\documentclass[%
 aip,cha,
 amsmath,amssymb,
 preprint,%
]{revtex4-1}

\usepackage{graphicx}% Include figure files
\usepackage{dcolumn}% Align table columns on decimal point
\usepackage{bm}% bold math
\usepackage{mathrsfs,subcaption}
\usepackage{color}

\def\rd{\text{d}}
\def\Zh{\widehat Z}
\def\eps{{\varepsilon}}

\def\fr{\mbox{$\frac{1}{2}$}}
\def\lth{\langle\!\langle}
\def\rth{\rangle\!\rangle}

\def\D{\mbox{${\rm D}$}}

\def\poa{\mbox{$\partial_{\omega_1}$}}
\def\pob{\mbox{$\partial_{\omega_2}$}}
\def\bw{{\bm \omega}}

\def\bk{{\bf k}}
\def\bO{{\bm \Omega}}
\def\be{{\bm \zeta}}
\def\bth{{\bm \theta}}

\def\bgam{{\bm \gamma}}

\def\bdel{{\bm \delta}}

\begin{document}

%\preprint{AIP/123-QED}

\title[Multiphase Modulation and the mKdV Equation]{The Modulation of Multiple Phases Leading to the Modified KdV Equation}% Force line breaks with \\

\author{D.J. Ratliff}
\email{d.j.ratliff@lboro.ac.uk}
\affiliation{ 
Department of Mathematical Sciences, Loughborough University, Loughborough, Leicestershire, United Kingdom, LE11 3TU%\\This line break forced with \textbackslash\textbackslash
}%

\date{\today}% It is always \today, today,

\begin{abstract}
This paper seeks to derive the modified KdV (mKdV) equation using a novel approach from systems generated from abstract Lagrangians possessing a two-parameter symmetry group. The method utilises a modified modulation approach, which results in the mKdV emerging with coefficients related to the conservation laws possessed by the original Lagrangian system. Alongside this, an adaptation of the method of Kuramoto is developed, providing a simpler mechanism to determine the coefficients of the nonlinear term. The theory is illustrated using two examples of physical interest, one in stratified hydrodynamics and another using a coupled Nonlinear Schr\"odinger model, to illustrate how the criterion for the mKdV equation to emerge may be assessed and its coefficients generated. 

\end{abstract}

%Use showkeys class option if keyword
\maketitle

 \begin{quotation}
Interacting nonlinear waves of two or more phases are a rich source of instability, which lead to the development of defects which then evolve over time to form further coherent structures, such as solitary pulses or nonlinear periodic forms. We present here one way in which the evolution of these defects can be modelled, by using the method of modulation to derive nonlinear partial differential equations which govern their evolution. In particular, we extend previous studies in this context to show that one may obtain a modified KdV (mKdV) equation, whose coefficients come from derivatives of the conservation of wave action associated with the original wavetrain. To help illustrate how this approach can be applied in practice, we study two physically relevant systems, a stratified shallow water system and a coupled Nonlinear Schr\"odinger model, in order to show how the conditions for the mKdV equation to emerge can be found and the relevant coefficients calculated.
 \end{quotation}

\section{Introduction}
The discussion of this paper centres around the modified Korteweg - de Vries (mKdV) equation, defined as 
\begin{equation}\label{mKdV}
a_0 q_T+a_1q^2q_X+a_2q_{XXX}= 0\,,
\end{equation}
for some unknown function $q(X,T)$ and coefficients $a_i$. This equation arises as a nonlinear reduction across various systems of interest, such as in interfacial flows \cite{cc99,g99,gpt97}, plasma physics \cite{ko69,nt94,ki74,s73} and thin ocean jets \cite{cpr92,rp94}. Moreover, it possesses several interesting solution families such as solitary waves, rational solutions and breathers that make this equation desirable to study \cite{zzsz14}. The interest of this paper is not in the solutions to this system, however, it is in using a new approach to derive (\ref{mKdV}) from systems generated from a Lagrangian density with two symmetries. Moreover, another aim of the paper is to show a connection between the coefficients of the mKdV equation and the conservation laws that the original Lagrangian system possesses.

The approach used to obtain the mKdV in this paper will be phase modulation. The roots of this approach are based in the works of \textsc{Whitham} \cite{WLNLW,w65}, who for single phased wavetrains derived the celebrated Whitham equations.
The theory starts by considering the abstract Lagrangian
\[
\mathscr{L}(U,U_x,U_t) = \iint \mathcal{L}(U,U_x,U_t) \ dx\,dt\,,
\]
for state vector $U(x,t)\in \mathbb{R}^n$ and Lagrangian density $\mathcal{L}$. One then assumes a periodic wavetrain solution to the associated Euler-Lagrange equations of the form 
\[
U = \hat{U}(kx+\omega t) \equiv \hat{U}(\theta; k,\omega), \quad \hat{U}(\theta+2 \pi) = \hat{U}(\theta)\,,
\]
exists, and so the Lagrangian averaged over one period becomes
\[
\mathscr{L}(k,\omega) = \iint \mathcal{L}(U,k U_\theta,\omega U_\theta) \ d \theta\,.
\]
The strategy is to then assume the phase, wavenumber and frequency are all slowly varying functions, so that $k = \theta_X$ and $\omega = \theta_T$ for slow variables $X = \eps x, \, T = \eps t$. Then by taking variations of the averaged Lagrangian with respect to $\theta$, one obtains
\[
(\mathscr{L}_\omega(k,\omega))_T+(\mathscr{L}_k(k,\omega))_X \equiv A(k,\omega)_T+B(k,\omega)_X = 0\,.
\]
This coupled with the consistency condition $k_T = \omega_X$, form the Whitham modulation equations, a set of dispersionless nonlinear partial differential equations (PDEs). It transpires that the functions $A$ and $B$ are the components of the conservation of wave action the Lagrangian possesses evaluated on the wavetrain solution $\hat{U}$. The Whitham modulation equations have since been obtained from a broader class of solutions known as \emph{relative equilibria}, which are solutions that are steady relative to the orbit of some group action \cite{BSNLW}. This generalises the modulation of periodic waves, allowing one to consider a larger number of problems. It is for this reason that the modulation of relative equilibria forms the focus of the paper.

These ideas can be extended to wavetrains with multiple phases to recover similar results \cite{ab70,r17a}. One may repeat the above procedure, but instead consider the two-phased doubly periodic wavetrain (and in general, two-phased relative equilibrium)
\[
\begin{split}
&U =\hat{U}(\bth; \bk, \bw)\,, \quad\hat{U}(\theta_1+2 \pi,\theta_2) = \hat{U}(\bth) = \hat{U}(\theta_1,\theta_2+2 \pi)\,,\\
& \bth = 
\begin{pmatrix}
\theta_1\\
\theta_2
\end{pmatrix} = 
\begin{pmatrix}
k_1x+\omega_1t\\
k_2 x+\omega_2 t
\end{pmatrix}\,, \quad
 \bk = 
\begin{pmatrix}
k_1\\
k_2
\end{pmatrix}\,, \ \bw
\begin{pmatrix}
\omega_1\\
\omega_2
\end{pmatrix}\,.
 \end{split}
\]
By assuming each phase is slowly varying again, so that $\bk = \bth_X$ and $\bw = \bth_T$, variations of the $\bth$-averaged Lagrangian lead to the vector Whitham modulation equations:
\[
A(\bk,\bw)_T+B(\bk,\bw)_X = {\bf 0}\,, \quad \bk_T = \bw_X\,.
\]
In this case, $A$ and $B$ are vector valued, and their components form the conservation of wave action associated with each phase.

An interesting avenue of research has focussed on the case where the Whitham equations are degenerate. This corresponds to the emergence of a zero characteristic in its linearisation. For the scalar Whitham equation, this zero characteristic emerges at points where $B_k(k_0,\omega_0) \equiv \mathscr{B}_k = 0$ for fixed wavenumber and frequency $k_0,\,\omega_0$. At such points it has been shown that a more general modulation approach is required. This has been developed by \textsc{Bridges} \cite{b13}, and was in part inspired by the work of \textsc{Doelman et. al.} \cite{DSSS}. The idea is to construct a \emph{modulation ansatz}, which takes the relative equilibrium solution and perturbs each of its independent variables:
\begin{multline}\notag
U = \hat{U}\big(\theta+\eps \phi(X,T), k +\eps^2 q(X,T), \omega+\eps^4\Omega(X,T)\big)
+\eps^3W(\theta,X,T)\,,
\end{multline}
where $q = \phi_X,\,\Omega = \phi_T$ and $W$ is a remainder term, which is required since $\hat{U}$ is no longer an exact solution. The slow variables are rescaled as $X = \eps x,\, T =\eps^3 t$ in light of the zero characteristic. 
This guess at a solution is then substituted into the Euler-Lagrange equations associated with the abstract Lagrangian. By undertaking the resulting asymptotic analysis, dispersion arises from the modulation and leads to the emergence of the Korteweg- de Vries (KdV) equation:
\begin{equation}\label{KdVEqtn}
(\mathscr{A}_k+\mathscr{B}_\omega)q_T+\mathscr{B}_{kk}qq_X+\mathscr{K}q_{XXX} = 0\,.
\end{equation}
It is apparent that the majority of the coefficients are related to the conservation laws $A$ and $B$, and the dispersive coefficient $\mathscr{K}$ can be obtained from a Jordan chain analysis. 

The benefits of the above modulation approach are two-fold. Firstly, since the Lagrangian considered is abstract, the calculations for the reduction need only be done once in order to apply to any Lagrangian where the required criterion can be met. This means that the results that emerge have the potential to be widely applicable across many areas of physics. Secondly, the majority of the coefficients in the nonlinear PDE obtained from the modulation approach are related to derivatives of the conservation laws. These can be determined \emph{a priori} to the modulation, and the calculation of the necessary derivatives is typically elementary. This is the primary reason that a version of this method is adopted within this paper, so that these benefits may also apply to the results presented within this work.

Along with the degeneracy of the Whitham equations, the KdV equation (\ref{KdVEqtn}) may also degenerate, occuring when one of (or combinations of) its coefficients vanish. A similar set of generalised modulation approaches show that in such scenarios other well known nonlinear PDEs emerge from the analysis \cite{rb16,r17,rb17}. The last of these shows that in cases where $\mathscr{B}_k = \mathscr{B}_{kk} = 0$ the analysis admits the mKdV equation in the form
\begin{equation}\label{mKdVSingle}
(\mathscr{A}_k+\mathscr{B}_\omega)q_T+\frac{1}{2}\mathscr{B}_{kkk}q^2q_X+\mathscr{K}q_{XXX} = 0\,.
\end{equation}
Once again the conservation laws are related to the majority of its coefficients, albeit a higher derivative is now present as the coefficient of the new nonlinear term. The principle aim of this paper will be to generalise this single phase result to the case of two phases in a way that can also be extended to arbitrarily many.

The modulation of multiple phases in the presence of zero characteristics has very recently been developed. In the case of the vector Whitham equations, the emergence of a zero characteristic may be shown to occur precisely when
\begin{equation}\label{BkCond}
{\rm det}\big[\D_\bk B(\bk_0,\bw_0)\big] \equiv {\rm det}\big[\D_\bk {\bf B}\big] = 0\,,
\end{equation}
where $\D$ denotes the directional derivative with respect to the subscripted argument and the bold notation denotes evaluation at the constant vectors $\bk_0,\bw_0$. This generalises the single phase condition naturally and allows one to define the eigenvector $\be$ with the property
\[
\D_\bk{\bf B}\be = {\bf 0}\,.
\]
To abridge the single phase approach one constructs the modulation ansatz
\begin{multline}
U = \hat{U}\big(\bth+\eps \be \phi(X,T),\bk+\eps^2 \be U(X,T),\bw+\eps^4 \be \Omega(X,T)\big)
+\eps^3W(\bth,X,T)\,,
\end{multline}
where $U = \phi_X,\,\Omega = \phi_T$ and the slow variables are again scaled as $X = \eps x,\, T = \eps^3 t$.  Upon substitution of this expression into the Euler-Lagrange equations, one is able to show that when the above condition holds a vector KdV-like equation emerges from the analysis:
\[
\begin{split}
(\D_\bk{\bf A}_\bk+\D_\bw{\bf B})\be& U_T+\D_\bk^2{\bf B}(\be,\be)UU_X+{\bf K}U_{XXX}+\D_\bk{\bf B}{\bm \alpha}_{XX} = {\bf 0}
\end{split}
\]
for unknown function $U(X,T)$ and arbitrary vector-valued function ${\bm \alpha}(X,T)$ required to ensure the analysis results in nontrivial $U$ \cite{rb16a}. This can be turned into the scalar KdV equation by multiplying on the left by $\be$, which removes the ${\bm \alpha}$ term and gives the KdV equation
\begin{multline}\label{KdVMulti}
\be^T(\D_\bk{\bf A}_\bk+\D_\bw{\bf B})\be U_T+\be^T\D_\bk^2{\bf B}(\be,\be)UU_X\\+\be^T{\bf K}U_{XXX} = 0\,.
\end{multline}
Once again, there is a connection between the conservation laws evaluated along the solution and the coefficients of the resulting KdV. This paper is concerned with one of the cases that lead to the above KdV being degenerate, which will be when the coefficient of the nonlinear term vanishes and so
\begin{equation}\label{BkkCond}
\be \D_\bk^2{\bf B}(\be,\be) = 0\,.
\end{equation}
The results of the studies for single phase modulation would suggest that the analysis in this case would lead to the mKdV, and the main result of this paper confirms this. It will be shown that the modulation approach in light of the conditions (\ref{BkCond}) and (\ref{BkkCond}) holding simultaneously leads to the emergence of an mKdV equation of the form
\begin{multline}\label{mKdVUniIntro}
\be^T(\D_\bk{\bf A}_\bk+\D_\bw{\bf B})\be V_T+\frac{1}{2}\be^T\big(\D_\bk^2{\bf B}(\be,\be,\be)-3\D_\bk^2{\bf B}(\be,\bdel)\big)V^2V_X+\be^T{\bf K}V_{XXX} = 0\,,
\end{multline}
for unknown function $V(X,T)$ and the vector $\bdel$ satisfies
\[
\D_\bk{\bf B}\bdel = \D_\bk^2{\bf B}(\be,\be)\,.
\]
The similarities between (\ref{mKdVUniIntro}) and (\ref{mKdVSingle}) are quite clear, although the generalisation is not entirely trivial due to the presence of the $\bdel$ term. The modulation analysis presented in this paper will emphasise the role of this vector and how it arises in the theory.

In order to justify the new form of the nonlinear coefficient, and another key contribution of this paper, we develop a method to determine the nonlinear coefficient of the resulting modulation equation without the need to undertake the modulation analysis. This is achieved by adapting the method of Kuramoto used in the modulation of single phase wavetrains \cite{k84} to multiple phases. The essence of the method is that the coefficients of the nonlinearity arise from Taylor expansions of the Whitham equations, and the idea for the multiphase case is no different. This extension is somewhat natural, with one instead dealing with the derivatives of tensors instead of scalars, meaning that the results are somewhat identical. The calculations involved are somewhat easier than those resulting from the modulation analysis, however the two are shown to be in agreement. Overall, this development provides an easier avenue to generate the coefficients of the nonlinear terms obtained from the modulation approach. This extended method of Kuramoto, although developed to validate the mKdV derived here, can be used to obtain coefficients across several other modulation analyses.

To demonstrate how the result of this paper may be applied, we illustrate two examples of how the mKdV equation may arise from two physically relevant systems. The first is by using a stratified shallow water system, where the mKdV equation is shown to emerge from flows of finite speed providing the relevant criteria are met. This is also a step forward from the literature, where the mKdV is typically derived for flows of zero velocity \cite{dr78,kb81,g99}. The other example considered in this paper is a pair of coupled Nonlinear Schr\"odinger (NLS) equations, where it will be shown that the mKdV equation may be obtained via the modulation of plane waves. This is the first such reduction from the coupled NLS system to the scalar mKdV that the author is aware of, and so the theory presented in this paper leads to the emergence of the mKdV equation in new contexts.

The structure of this paper is as follows. In \S\ref{sec:AS} the relevant abstract setup for the theory is developed. Within this, properties of the relative equilibrium solution, the structure of the conservation laws and the relevant Jordan chain theory are discussed. This is followed by the extension of the method of Kuramoto to multiphase modulation in \S\ref{sec:kura}, showing how the coefficients of nonlinear terms can be obtained by considering Taylor expansions of the fully nonlinear Whitham equations. The modulation analysis leading to the mKdV equation is presented in \S\ref{sec:ModAnal}, demonstrating how the two conditions (\ref{BkCond}), (\ref{BkkCond}) result in the equation (\ref{mKdVUniIntro}) emerging. Examples of how the theory applies to problems of interest are given in \S\ref{sec:SWH} and \ref{sec:NLS}, demonstrating how the mKdV equation arises from both stratified shallow water hydrodynamics and a coupled Nonlinear Schr\"odinger model. Concluding remarks are given at the end of the paper.

\section{Abstract setup}\label{sec:AS}
The starting point for the theory of this paper is the class of problems generated by a Lagrangian density. In particular, we make the assumption that this density is in multisymplectic form. The process of transforming a Lagrangian into multisymplectic form is essentially a sequence of Legendre transformations, which are documented in detail in another work \cite{rb16}, and so this is not recounted here. Instead, we state that the multisymplectic Lagrangian takes the form
\begin{equation}\label{MSFL}
\mathscr{L} = \iint \bigg(\frac{1}{2}\langle Z, {\bf M}Z_t\rangle+\frac{1}{2}\langle Z, {\bf J}Z_x\rangle - S(Z)\bigg)\ dx\,dt\,,
\end{equation}
for state vector $Z\in \mathbb{R}^n$, $\langle \cdot,\cdot \rangle$ denotes the standard inner product on $\mathbb{R}^n$, ${\bf M}$, ${\bf J}$ are constant skew-symmetric matrices and $S$ denotes some Hamiltonian function which is generated through the Legendre transformations. The Euler-Lagrange equations for the system are obtained by taking the first variation of the Lagrangian density, which for the multisymplectic Lagrangian above gives
\begin{equation}\label{MSFELE}
{\bf M}Z_t+{\bf J}Z_x = \nabla S(Z)\,.
\end{equation}
This system will be one of the main constructs discussed in this paper, as it will be solutions to this equation that are modulated and the mKdV will be obtained as a reduction to this system.

The methodology of this paper proceeds under the assumption that the system (\ref{MSFELE}) possesses a two phased relative equilibrium solution. Relative equilibria are solutions associated with a continuous symmetry which move along the orbit of the group. Such solutions can be thought of as the generalisation of wavetrains with two phases, which themselves are solutions associated with the invariance of phase translations. These solutions are of the form
\begin{multline}
Z(x,t) = \Zh(\theta_1,\theta_2,k_1,k_2,\omega_1,\omega_2) \equiv \Zh(\bth,\bk,\bw)\,, \\
\bth = \begin{pmatrix}
\theta_1\\
\theta_2
\end{pmatrix}\,,  \bk= 
\begin{pmatrix}
k_1\\
k_2
\end{pmatrix}\,, \ \bw = 
\begin{pmatrix}
\omega_1\\
\omega_2
\end{pmatrix}\,.
\end{multline}
The wavenumbers $k_i$ and frequencies $\omega_i$ are taken to be constant in these solutions. Substitution of this expression into (\ref{MSFELE}) generates the PDE
\begin{equation}\label{BasicState}
\sum_{i=1}^2\big(\omega_i{\bf M}+k_i{\bf J}\big)\Zh_{\theta_i} = \nabla S(\Zh)\,.
\end{equation}

The linearisation of the above PDE arises frequently within the modulation analysis, which allows one to define the associated linear operator ${\bf L}$ as
\[
{\bf L}V = \D^2S(\Zh)-\sum_{i=1}^2\big(\omega_i{\bf M}+k_i{\bf J}\big)V_{\theta_i}\,.
\]
In particular, the operator ${\bf L}$ is self adjoint under the $\bth$-averaging inner product
\[
\lth U,\,V\rth = \frac{1}{4\pi^2}\int_0^{2 \pi}\int_0^{2 \pi} \langle U,\,V\rangle \ d \theta_1\,d \theta_2 \quad \forall \ U,\,V \in \mathbb{R}^n\,.
\]
For symmetries that are affine (such as the first example of this paper) the averaging is dropped and this becomes the standard inner product on $\mathbb{R}^n$. By differentiating (\ref{BasicState}) with respect to each of the parameters $\theta_i,\,k_i$ and $\omega_i$, one is able to obtain the following results:
\begin{subequations}
\begin{align}
{\bf L} \Zh _{\theta _i} &= 0, \label{ThetaDeriv}\\
{\bf L} \Zh _{k_i} &= {\bf J} \Zh _{\theta _i}, \label{kDeriv}\\
{\bf L} \Zh _{\omega _i} &= {\bf M} \Zh _{\theta _i}. \label{OmegaDeriv}
\end{align}
\end{subequations}
The first of these equations highlights that each of the $\Zh_{\theta_i}$ lie within the kernal of ${\bf L}$. An assumption made in this paper is that this kernel is no larger. This means that the solvability requirement for inhomogenous problems takes the form
\begin{equation}\label{solvability}
{\bf L}F = G \quad \mbox{is only solvable when} \quad \lth \Zh_{\theta_i},G\rth = 0\,, \quad i = 1,\,2\,,
\end{equation}
The remaining two equations, (\ref{kDeriv}) and (\ref{OmegaDeriv}), highlight that Jordan chains arise involving ${\bf L}$. Of these, only the one involving the matrix ${\bf J}$ is important in this paper and the theory for such chains is reviewed in \S\ref{sec:jct}. This chain will ultimately be the mechanism for which dispersion enters the phase dynamics.

\subsection{Symmetries and conservation laws}
One benefit of putting the Lagrangian in the form (\ref{MSFL}) is that an explicit connection between the system's conservation laws and the structure of the Euler-Lagrange equations can be made. This is through the symplectic structures ${\bf M}$ and ${\bf J}$, which appear in both the Euler-Lagrange equations and the conservation laws. By appealing to Noether theory for multisymplectic Lagrangians \cite{bhj10,h05} in the case of two symmetries, the conservation laws may be found as
\[
A(x,t) = \frac{1}{2}
\begin{pmatrix}
\lth Z,{\bf M}Z_{s_1}\rth\\
\lth Z,{\bf M}Z_{s_2}\rth
\end{pmatrix}\,, \quad B(x,t) = 
 \frac{1}{2}
\begin{pmatrix}
\lth Z,{\bf J}Z_{s_1}\rth\\
\lth Z,{\bf J}Z_{s_2}\rth
\end{pmatrix}
\]
where $s_i$ parameterise each of the symmetries associated with the solution. In the case of relative equilibria considered in this paper, we simply have that $s_i = \theta_i$. The affine case is almost identical but without the factors of $\frac{1}{2}$. One is able to evaluate these along the solution $\Zh$ to obtain the vectors
\[
\begin{split}
{\bf A}(\bk,\bw) =& \begin{pmatrix} \mathscr{A}_1 \\
\mathscr{A}_2\end{pmatrix}
: =\frac{1}{2}
\begin{pmatrix} \lth{\bf M}\Zh_{\theta_1},\Zh\rth \\
 \lth{\bf M}\Zh_{\theta_2},\Zh\rth \end{pmatrix}\,,
 \\
{\bf B}(\bk,\bw) =& \begin{pmatrix} \mathscr{B}_1 \\
\mathscr{B}_2\end{pmatrix} := \frac{1}{2}
\begin{pmatrix} \lth{\bf J}\Zh_{\theta_1},\Zh\rth \\
 \lth{\bf J}\Zh_{\theta_2},\Zh\rth \end{pmatrix}\,.
 \end{split}
\]
 These expressions in the periodic case can also be obtained through the $k$ and $\omega$ derivatives of the Lagrangian
(\ref{MSFL}) averaged over the two-phase solution:
\begin{multline}\notag
\mathscr{L}(\bk,\bw) = \frac{1}{4\pi^2}\int_0^{2\pi}\!\!\!\!\int_0^{2\pi}
\Bigg[\frac{1}{2}\sum_{j=1}^2\left[
 \langle \Zh, \omega_j{\bf M}\Zh_{\theta_j} +k_j {\bf J} \Zh_{\theta_j}\rangle \right]
-S(\Zh) \Bigg]\,\rd\theta_1\rd\theta_2\,.
\end{multline}

By the definitions of these conservation laws, one is able to obtain the following tensors of derivatives:
\begin{equation}\notag
\begin{split}
{\rm D}_{\bf k}{\bf A} =& \begin{pmatrix}
\partial_{k_1}\mathscr{A}_1&\partial_{k_2}\mathscr{A}_1\\
\partial_{k_1}\mathscr{A}_2&\partial_{k_2}\mathscr{A}_2
\end{pmatrix} = {\rm D}_{\bm \omega}{\bf B}^T,\\[2mm]
\D_\bw {\bf A} =& 
\begin{pmatrix}
\poa \mathscr{A}_1&\pob \mathscr{A}_1\\
\poa\mathscr{A}_2&\pob \mathscr{A}_2
\end{pmatrix}, \quad
{\rm D}_{\bf k}{\bf B} = \begin{pmatrix}
\partial_{k_1}\mathscr{B}_1&\partial_{k_2}\mathscr{B}_1\\
\partial_{k_1}\mathscr{B}_2&\partial_{k_2}\mathscr{B}_2
\end{pmatrix},\\[2mm]
{\rm D}^2_{\bf k}{\bf B} =& \left(\begin{array}{cc}
\partial_{k_1 k_1}\mathscr{B}_1&\partial_{k_2 k_1}\mathscr{B}_1\\
\partial_{k_1 k_1}\mathscr{B}_2&\partial_{k_2 k_1}\mathscr{B}_2
\end{array}\vline
\begin{array}{cc}
\partial_{k_1 k_2}\mathscr{B}_1&\partial_{k_2 k_2}\mathscr{B}_1\\
\partial_{k_1k_2}\mathscr{B}_2&\partial_{k_2k_2}\mathscr{B}_2
\end{array}\right)\,.
\end{split}
\end{equation}
The entries of these tensors are related to the solution $\Zh$ by
\begin{subequations}
\begin{align}
\partial _{k_j} \mathscr{A}_i &= \lth {\bf M} \Zh _{\theta _i} , \Zh _{k_j} \rth , \label{Ak} \\
\partial _{\omega_j} \mathscr{A}_i &= \lth {\bf M} \Zh _{\theta _i} , \Zh _{\omega_j} \rth , \label{Aw}\\
\partial _{k_j} \mathscr{B}_i &= \lth {\bf J} \Zh _{\theta _i} , \Zh _{k_j} \rth , \label{Bk} \\
\partial _{k_j k_m} \mathscr{B} _i & = \lth {\bf J} \Zh_{\theta _i k_m} , \Zh _{k_j} \rth + \lth {\bf J} \Zh _{\theta _i} , \Zh _{k_j k_m} \rth.\label{Bkk}
\end{align}
\end{subequations}
We note that
\begin{multline}\label{BDerivRelation}
\partial _{k_i} \mathscr{B} _j = \lth {\bf J} \Zh _{\theta _j} , \Zh _{k_i} \rth = \lth {\bf L} \Zh _{k _j} , \Zh _{k_i} \rth = \lth \Zh _{k _j} , {\bf L} \Zh _{k_i} \rth = \lth \Zh _{k _j} , {\bf J} \Zh _{\theta _i} \rth = \partial _{k_j} \mathscr{B}_i
\end{multline}
along with
\begin{equation}\notag
\partial _{k_j}\mathscr{A}_i = \lth {\bf M}\Zh_{\theta _i} , \Zh _{k_i} \rth = \lth \Zh _{\omega _i} , {\bf J} \Zh _{k _j}\rth = \partial _{\omega _i}\mathscr{B}_j\,.
\end{equation}

The notion of criticality plays a fundamental role in the modulation approach,
 as it is the mechanism that leads to the emergence of nonlinear dynamics.
In the context of this paper, we define that a conservation law is critical 
if it's Jacobian with respect to either $\bk$ or $\bw$ has a zero determinant.
Criticality in this sense then holds along surfaces in $(\bk,\,\bw$)-space, and the modulation equations that emerge from this theory are valid along such curves or sufficiently close to them. The primary criticality this paper is concerned with arises when
\begin{equation}\label{DetBZero}
{\rm det}\big[\D_\bk{\bf B}\big] = 0\,,
\end{equation}
 which corresponds to the emergence of a zero characteristic from the Whitham equations obtained from the Lagrangian \cite{rb16,r17}. It also facilitates the definition of the eigenvector associated with this zero eigenvalue, denoted as $\be$, so that
\begin{equation}
\D_\bk{\bf B}\be = {\bf 0}\,.
\end{equation} 
Throughout the paper the zero eigenvalue is assumed to be simple, so that there is only one such kernel element, although the theory may be abridged when this is not true. There is a link between the condition (\ref{DetBZero}) and the emergence of dispersion from the modulation approach, which is discussed in \S\ref{sec:jct}. Interestingly, this condition also arises across the literature as a stability boundary \cite{blp01,b54,l90}, and so the emergence of nonlinear PDEs has an interesting connection to the stability of the system.
 
 This paper extends the notion of criticality further by considering the case where the second directional derivative of ${\bf B}$ vanishes in the direction of $\be$, meaning that
 \begin{equation}\label{BkkZeroCons}
 \be^T\D_\bk{\bf B}(\be,\be) = 0\,.
 \end{equation}
This is precisely when the nonlinear term in the KdV given in (\ref{KdVMulti}) vanishes, which would imply that the modulation approach needs to be altered in such cases. This rescaling is undertaken in \S\ref{sec:ModAnal}. The condition (\ref{BkkZeroCons}) also arises as the condition that the system
\[
\D_\bk {\bf B}\bdel = \D_\bk{\bf B}(\be,\be)\,,
\]
is solvable, since $\be$ lies in the kernel of $\D_\bk {\bf B}$. This will be how the additional vector $\bdel$ enters into the modulation analysis leading to the additional term in (\ref{mKdVUniIntro}). The precise details of this will be revisited in \S\ref{sec:ModAnal}.

\subsection{Review of the jordan chain theory for multiple phases}\label{sec:jct}
The other construct arising from the modulation approach is a Jordan chain analysis, as suggested by the results (\ref{kDeriv}) and (\ref{OmegaDeriv}). We review the relevant Jordan chain theory generated by the former, since this will be the mechanism that leads to the emergence of dispersion from the phase dynamics.

We can see from (\ref{ThetaDeriv}) and (\ref{kDeriv}) that we begin to form two Jordan chains with the structure
\begin{equation}\notag
{\bf L}\xi _1 = 0, \quad {\bf L}\xi _i = {\bf J} \xi _{i-1}, \ i>1.
\end{equation}
The two chains are started with the $\theta$ derivatives and are followed by the respective $k$ derivative. We denote these in the following way:
\[
\xi_1 = \Zh_{\theta_1}\,, \quad \xi_2 = \Zh_{k_1}\,, \quad \xi_3 = \Zh_{\theta_2}\,, \quad \xi_4 = \Zh_{k_2}\,,
\]
so that the first two form the first chain, and the latter two the second. However, these chains will coalesce to allow the modulation analysis to continue. Consider the equation
\begin{equation}\label{xi5eqtn}
{\bf L}\xi_5 = \sum_{i=1}^2 \zeta_i{\bf J}\Zh_{k_i}\,.
\end{equation}
Assessing the solvability of the above generates the system
\begin{multline}\label{EigenProblem}
\begin{pmatrix}
\lth \Zh_{\theta_1},{\bf J}\Zh_{k_1}\rth&\lth \Zh_{\theta_1},{\bf J}\Zh_{k_2}\rth\\
\lth \Zh_{\theta_2},{\bf J}\Zh_{k_1}\rth&\lth \Zh_{\theta_2},{\bf J}\Zh_{k_2}\rth
\end{pmatrix}
\begin{pmatrix}
\zeta_1\\
\zeta_2
\end{pmatrix}\\ \equiv -
\begin{pmatrix}
\partial_{k_1}\mathscr{B}_1&\partial_{k_2}\mathscr{B}_1\\
\partial_{k_1}\mathscr{B}_2&\partial_{k_2}\mathscr{B}_2
\end{pmatrix}
\be \equiv-\D_\bk{\bf B}\be = {\bf 0}\,.
\end{multline}
Therefore, if $\D_\bk {\bf B}$ has a zero eigenvalue with eigenvector $\be$, then the above system (\ref{xi5eqtn}) is solvable. This also requires that (\ref{DetBZero}) holds so that the matrix possesses a zero eigenvalue. In such cases, one is able to define
\begin{equation}\notag
{\bf L} \xi_5 = \sum_{i=1}^2\zeta_i{\bf J} \Zh _{k_i}.
\end{equation}
The system (\ref{EigenProblem}) has another eigenvalue given by the trace of $\D_\bk{\bf B}$, which results in the eigenvalue problem
\begin{equation}\label{NonZeroEigen}
\D_\bk {\bf B}
\left(
\begin{array}{c}
\zeta_2\\
-\zeta_1
\end{array}
\right) = (\partial _{k_1} \mathscr{B}_1+ \partial _{k_2} \mathscr{B}_2)
\left(
\begin{array}{c}
\zeta_2\\
-\zeta_1
\end{array}
\right).
\end{equation}
A consequence of the above is that the equation
\begin{equation}\notag
{\bf L}F = \zeta_2{\bf J}\Zh _{k_1} - \zeta_1 {\bf J} \Zh _{k _2}\,,
\end{equation}
is no longer solvable, as the zero eigenvalue is assumed simple.

 Because the zero eigenvalue of ${\bf L}$ is even, the existence of $\xi_5$ guarantees the existence of $\xi_6$ with
\begin{equation}\label{xi6}
{\bf L}\xi_6 = {\bf J}\xi_5\,.
\end{equation}
In particular, the fact that this system is solvable gives that
\begin{equation}\label{k-k-zero}
\begin{split}
0=\lth \Zh_{\theta_i},{\bf J}\xi_5\rth = -\sum_{j=1}^2\zeta_j\lth \Zh_{k_i},{\bf J}\Zh_{k_j}\rth\,,\\
\therefore \quad \lth \Zh_{k_i},{\bf J}\Zh_{k_j}\rth = 0\,, \quad i,\,j=1,\,2\,.
\end{split}
\end{equation}
and the assumption in this paper will be that this chain is no longer. As a consequence, the system
\[
{\bf L}\xi_7 = {\bf J}\xi_6\,,
\]
is not solvable, which means that
\begin{equation}\label{BoldKDefn}
{\bf K}  = 
\begin{pmatrix}
\mathscr{K}_1\\
\mathscr{K}_2
\end{pmatrix}
:= -
\begin{pmatrix}
\lth \Zh_{\theta_i},{\bf J}\xi_6\rth\\
\lth \Zh_{\theta_2},{\bf J}\xi_6\rth
\end{pmatrix} \neq {\bf 0}\,.
\end{equation}
This vector forms the dispersive component in the mKdV derived here, and will arise directly from the phase dynamics.

\section{Method of Kuramoto in multiphase modulation}\label{sec:kura}
The calculations arising from the modulation in the context of this paper will generate several involved calculations, and so the question arises as to how accurate these results might be. In order to confirm these, as well as present an alternate ad-hoc way for which these can be obtained, we abridge a technique for obtaining nonlinear coefficients from the modulation single phase wavetrains in order to use it for the analysis presented in this paper.

The method of Kuramoto provides a useful tool when discussing the coefficients of nonlinearities in phase dynamics \cite{hpf,k84}. The technique was originally developed in non-conservative single phase modulation, illustrating how one may deduce the coefficients of nonlinear terms in the reduced equation without requiring further modulation. Inspired by this technique, one is able to modify the approach slightly for the case of tensors. The principle remains the same, and these modifications are detailed below.

Consider the fully nonlinear multiphase Whitham equations:
\begin{equation}\label{NonWhit}
{\bf A}({\bm K},\bO)_T+{\bf B}({\bm K},\bO)_X = {\bf 0}\,,
\end{equation}
where  ${\bm K},\,\bO$ are the slowly varying wavenumber and frequency which are functions of $X,\,T$ whose scales at this stage are $X = \eps x,\, T = \eps t$. Consider now an expansion of the slowly varying wavenumber of the form
\[
{\bm K} = \bk+\sum_{i=1}^\infty\frac{\eps^n}{n!}U^n(X,T){\bm \chi}_n\,,
\]
for $\bk,\,{\bm \chi}_n$ some fixed vectors, $U$ some slowly varying function and $\eps \ll 1$. For simplicity we also fix $\bO$ as some constant vector. The scale of the slow variable $X$ will remain the same, but the scaling of $T$ will depend on the dispersive term present, which is not discussed here. Substituting the above into the function ${\bf B}$ in (\ref{NonWhit}) initially and Taylor expanding about $\eps = 0$ morphs the ${\bf B}$ term as
\begin{widetext}
\begin{multline}\notag
{\bf B}({\bm K},\bO)_X = \bigg(\eps\D_\bk {\bf B}{\bm \chi}_1U+\frac{1}{2}\eps^2U^2\big(\D_\bk{\bf B}{\bm \chi}_2+\D_\bk^2{\bf B}({\bm \chi}_1,{\bm \chi}_1)\big)+\frac{1}{6}\eps^3U^3\big(\D_\bk{\bf B}{\bm \chi}_3+3\D_\bk^2{\bf B}({\bm \chi}_1,{\bm \chi}_2)+\D_\bk^3{\bf B}({\bm \chi}_1,{\bm \chi}_1,{\bm \chi}_1)\big)\bigg)_X\\
+\mathcal{O}(\eps^4)\,.
\end{multline}
\end{widetext}
This expansion may be continued up to the order desired, depending on which terms are nonvanishing. The key idea of the method is to then consider which term in the largest set of brackets is the leading order term. For most cases this is the first term and the analysis becomes that of the linear Whitham equations \cite{r17a}, however in cases where
\[
{\rm det}\big[\D_\bk {\bf B}\big] = 0\,,
\]
then this term vanishes whenever ${\bm \chi}_1 = \be$, which will be used throughout the remainder of this discussion. 

The next term in this bracket then becomes important. If ${\bm \chi_2}$ is zero then what one would obtain is precisely the quadratic nonlinearity term obtained in the derivation of the KdV equation and the two-way Boussinesq equation via the modulation of multiple phases \cite{rb16,r17a}. There are however scenarios in which the projection of the nonlinearity vanishes along with the first term, meaning
\[
\be^T\D_\bk{\bf B}(\be,\be) = 0\,.
\]
This would imply that the quadratic nonlinearity of the scalar phase equations vanishes, and so a rescaling should occur to replace it. This condition implies that the system
\[
\D_\bk{\bf B}{\bm \chi}_2+\D_\bk^2{\bf B}(\be, \be) = 0\,,
\]
may be solved. For the purposes of this paper we do so by setting ${\bm \chi}_2 = -\bdel$ so that
\[
\D_\bk{\bf B}{\bm \delta} = \D_\bk^2{\bf B}(\be,\be)\,.
\]
The cubic terms are then the most dominant nonlinear term arising in the expansion, and so for ${\bm \chi}_3$ chosen to be zero (since this will be the highest order considered) the dominant term of the expansion is
\begin{equation}\label{cubickura}
{\bf B}({\bm K},\bO)_X = \frac{1}{2}\big(\D_\bk^3{\bf B}(\be,\be,\be)-3\D_\bk^2{\bf B}(\bdel,\be)\big)U^2U_X+\mathcal{O}(\eps^4)\,.
\end{equation}
The analysis presented in this paper will show that this is precisely the cubic term one obtains from the modulation approach. Thus, this method allows one to obtain the necessary term for the nonlinearity without having to undertake the modulation. This is expected to be true for higher nonlinearities (such as $U^3U_X$), and when $\bO$ is not fixed for nonlinearities of a mixed type (such as $UU_T$, $U_XU_{XT}$), although for this latter case one must also consider an expansion of ${\bf A}$.

\section{Summary of the modulation reduction}\label{sec:ModAnal}
We now present the detail of the modulation leading to (\ref{mKdV}). In order to achieve this, we will use the ansatz
\begin{multline}\label{Ansatz}
Z = \Zh\big(\bth+\be U(X,T;\eps)-\eps {\bm \Phi}(X,T,\eps),\bk+\eps \be U_X-\eps^2 {\bm \Phi}_X,\bw +\eps^3\be U_T+\eps^4 {\bm \Phi}_T\big)+\eps^2W(\bth+\be U,X,T;\eps)\,,
\end{multline}
with
\[
 {\bm \Phi} = {\bm \delta}P(X,T;\eps)+\eps{\bm \alpha}(X,T;\eps)\,.
\]
The function $P$ has the property
\[
P_X = \frac{1}{2}(U_X)^2\,,
\]
and $\be,\,{\bm \delta}$ satisfy the equations
\begin{equation}\label{deltadefn}
\D_\bk{\bf B}\be = 0\,, \quad \D_\bk^2{\bf B}(\be,\be) = \D_\bk{\bf B}{\bm \delta}\,.
\end{equation}
In order for the above systems to be solvable, we require the conditions that
\[
{\rm det}\big[\D_\bk{\bf B}\big] = 0\,, \quad \be^T\D_\bk^2{\bf B}(\be,\be) = 0\,.
\]
The function ${\bm \alpha}$ is considered arbitrary and used to motivate the final projection from a vector system to a scalar PDE. Only the leading order terms are needed of many of the slowly varying functions appearing above, aside from $W$, which is expanded as a series:
\[
W = \sum_{i=1}^\infty\eps^nW_n(\bth+\be U,X,T)\,,
\]
so that parts of the remainder term appear at each relevant order. We note its dependence on $\bth +\be U$ is due to the fact that the anstaz (\ref{Ansatz}) has leading order dependence on $U$ as well.

The advantage of incorporating these results in advance is three-fold - the first is that there is now only one important unknown function in the analysis, $U$, that will generate the emergent nonlinear PDE. The other, ${\bm \alpha}$ will be used to ensure the final matrix system is nontrivial and motivate the projection. Secondly is that by assuming the relevant conditions are met, the solvability condition at all orders apart from the last will happen automatically. Finally, the addition of these terms within the ansatz itself, rather than in $W$, lends itself to the cancellation of several unimportant terms due to the form of the ansatz and the multisymplectic structure. 

Below are the step to obtain the modified KdV in this setting, which emerges at fourth order in $\eps$. Although it emerges at such a low order within the analysis, one can see from the ansatz (\ref{Ansatz}) that solutions to the mKdV equation have leading order effect on the phase of the initial wavetrain, so the nonlinear effects are felt at leading order. We substitute the ansatz (\ref{Ansatz}) into the Euler-Lagrange equations (\ref{MSFELE}), Taylor expand around the $\eps = 0$ state and solve the system at each power of $\eps$. A summary of this is given below order by order.

\subsection{Leading order up to second order}
The leading order equation recovers the equation of the basic state (\ref{BasicState}). The next order gives that
\[
U_X\sum_{i=1}^2\zeta_i\big({\bf L}\Zh_{k_i}-{\bf J}\Zh_{\theta_i}\big) = 0\,,
\]
which is satisfied due to properties of the basic state. 

The next order, ignoring terms that cancel due to properties of $\Zh$, gives that
\[
{\bf L}W_0 = U_{XX}\sum_{i=1}^2\zeta_i{\bf J}\Zh_{k_i}\,.
\]
Applying the solvability condition (\ref{solvability}) gives that
\begin{equation}\label{BkSolv}
\begin{pmatrix}
\lth \Zh_{\theta_1},{\bf J}\Zh_{k_1}\rth & \lth \Zh_{\theta_1},{\bf J}\Zh_{k_2}\rth\\
\lth \Zh_{\theta_2},{\bf J}\Zh_{k_1}\rth & \lth \Zh_{\theta_2},{\bf J}\Zh_{k_2}\rth
\end{pmatrix}
\be \equiv - \D_\bk {\bf B} \be = {\bf 0}\,.
\end{equation}
As det$\big[\D_\bk{\bf B}\big] = 0$ is assumed, this holds by definition of $\be$, and so
\[
W_0 = U_{XX}\xi_5\,, \quad {\bf L}\xi_5 = \sum_{i=1}^2\zeta_i{\bf J}\Zh_{k_i}\,.
\]

\subsection{Third order}
The terms at third order, again ignoring those that cancel, gives
\[
\begin{split}
{\bf L}W_1 =& U_{XXX}{\bf J}\xi_5+U_XU_{XX}\sum_{i=1}^2\bigg[\zeta_i\big({\bf J}(\xi_5)_{\theta_i}\\
&-\D^3S(\Zh)(\Zh_{k_i}\xi_5)-\delta_i{\bf J}\Zh_{k_i}+\sum_{j=1}^2\zeta_j{\bf J}\Zh_{k_ik_j}\big)-\delta_i{\bf J}\Zh_{k_i}\bigg]\,.
\end{split}
\]
Appealing to solvability now, one can note that the $U_{XXX}$ term vanishes as the zero eigenvalue of ${\bf L}$ is even and so it is solvable. For the last term, we generate the system
\begin{widetext}
\begin{multline}
\begin{pmatrix}
\lth \Zh_{\theta_1},\sum_{i=1}^2\zeta_i\big({\bf J}(\xi_5)_{\theta_i}-\D^3S(\Zh)(\Zh_{k_i}\xi_5)+\sum_{j=1}^2\zeta_j{\bf J}\Zh_{k_ik_j}\big)\rth\\
\lth \Zh_{\theta_2},\sum_{i=1}^2\zeta_i\big({\bf J}(\xi_5)_{\theta_i}-\D^3S(\Zh)(\Zh_{k_i}\xi_5)+\sum_{j=1}^2\zeta_j{\bf J}\Zh_{k_ik_j}\big)\rth
\end{pmatrix} U_XU_{XX}\\
 - 
\begin{pmatrix}
\lth \Zh_{\theta_1},{\bf J}\Zh_{k_1}\rth & \lth \Zh_{\theta_1},{\bf J}\Zh_{k_2}\rth\\
\lth \Zh_{\theta_2},{\bf J}\Zh_{k_1}\rth & \lth \Zh_{\theta_2},{\bf J}\Zh_{k_2}\rth
\end{pmatrix}\bdel U_XU_{XX} = {\bf 0}\,, \\
 \implies \D_\bk{\bf B}{\bm \delta} = \D_\bk^2{\bf B}(\be,\be)\,.
\end{multline}
\end{widetext}
The full details of how the quadratic nonlinearity generates this term is given in appendix \ref{app:QuadCoeff}. This holds from the definition of ${\bm \delta}$, and so we may solve the problem at this order with
\[
W_1 = U_{XXX}\xi_6+U_XU_{XX}\kappa\,,
\]
with %$\xi_6$ defined in (\ref{xi6}) and
\[
{\bf L}\xi_6 = {\bf J}\xi_5
\]
as well as
\begin{multline}\notag
{\bf L}\kappa = \sum_{i=1}^2\bigg[\zeta_i\big({\bf J}(\xi_5)_{\theta_i}-\D^3S(\Zh)(\Zh_{k_i}\xi_5)+\sum_{j=1}^2\zeta_j{\bf J}\Zh_{k_ik_j}\big)-\delta_i{\bf J}\Zh_{k_i}\bigg]\,.
\end{multline}

\subsection{Fourth Order}\label{sec:ForuthOrd}
With the cancellation of many terms, the equation at this order reads
\begin{widetext}
\begin{equation}\label{fourthordmkdv}
\begin{split}
{\bf L}W_2 =&U_{XT}\sum_{i=1}^2\big({\bf M}\Zh_{k_i}+{\bf J}\Zh_{\omega_i}\big)+U_{XXXX}{\bf J}\xi_6+U_{XX}^2\big({\bf J}\kappa-\frac{1}{2}\D^3S(\Zh)(\xi_5,\xi_5)\big)\\
&+U_XU_{XXX}\big({\bf J}\kappa+\sum_{i=1}^2{\bf J}(\xi_6)_{\theta_i}-\D^3S(\Zh)(\Zh_{k_i},\xi_6)\big)+\sum_{i=1}^2(\alpha_i)_{XX}{\bf J}\Zh_{k_i}\\
&+U_{X}^2U_{XX}\sum_{i=1}^2\Bigg[-\frac{1}{2}\delta_i\big({\bf J}(\xi_5)_{\theta_i}-\D^2S(\Zh)(\xi_5,\Zh_{k_j})\big)+\zeta_i\bigg[{\bf J}(\kappa)_{\theta_i}-\D^3S(\Zh)(\kappa,\Zh_{k_i})\\
&-\sum_{j=1}^2\bigg(\frac{3}{2}\delta_j{\bf J}\Zh_{k_ik_j}+\frac{1}{2}\zeta_j\D^3S(\Zh)(\xi_5,\Zh_{k_ik_j})+\frac{1}{2}\D^4S(\Zh)(\xi_5,\Zh_{k_i},\Zh_{k_j})-\frac{1}{2}\sum_{m=1}^2\zeta_m\Zh_{k_i k_j k_m}\bigg)\bigg]\Bigg]\,.
\end{split}
\end{equation}
\end{widetext}
The idea is to now appeal to solvability and determine the tensors on each of the terms appearing in the above. For the $U_{XT}$ term, this generates
\[
\lth \Zh_{\theta_i},{\bf M}\Zh_{k_j}+{\bf J}\Zh_{\omega_j}\rth = -\partial_{k_j}\mathscr{A}_i-\partial_{\omega_j}\mathscr{B}_i\,,
\]
and for the $U_{XXXX}$ term, by definition,
\[
\lth \Zh_{\theta_i},{\bf J}\xi_6\rth = -\mathscr{K}_i\,.
\]
The terms involving $\alpha_i$ give
\[
\lth \Zh_{\theta_i},{\bf J}\Zh_{k_j}\rth = -\partial_{k_j}\mathscr{B}_i\,,
\]
as was seen in the computation undertaken in (\ref{BkSolv}). This completes the computation of the coefficients of the linear terms.
The quadratic nonlinearities at this order do not appear in the final PDE, since one is able to show that their coefficients are zero. This is expected due to their dissipative nature, and the relevant calculations for this can be found in appendix \ref{app:Vanishquad}. 

At this stage the equation governing solvability reads
\begin{multline}\notag
-(\D_\bk {\bf A}+\D_\bw {\bf B})\be U_{XT}+{\bf E}U_X^2U_{XX}\\- {\bf K}U_{XXXX}-\D_\bk {\bf B}{\bm \alpha}_{XX} = {\bf 0}\,.
\end{multline}
All that remains is to determine the coefficient of the $U_X^2U_{XX}$ term, the cubic nonlinearity. This calculation is considerably involved, but is undertaken in appendix \ref{app:CubicCoeff} and gives that
\[
{\bf E} = \frac{1}{2}\big(3\D_\bk^2{\bf B}(\bdel,\be)-\D_\bk^3{\bf B}(\be,\be,\be)\big)\,.
\]
This is in line with the coefficient generated by extending the method of Kuramoto in (\ref{cubickura}). 
With all terms accounted for, the final vector equation is:
\begin{multline}
(\D_\bk{\bf A}+\D_\bw{\bf B})\be U_{XT}+\frac{1}{2}\bigg(\D_\bk^3{\bf B}(\be,\be,\be)-3\D_\bk^2{\bf B}(\be,{\bm\delta})\bigg)U_X^2U_{XX}+{\bf K}U_{XXXX}+\D_\bk{\bf B}{\bm \alpha}_{XX} = 0\,.
\end{multline}
Using $\be$ to project this equation and the introduction of $V = U_X$ gives the scalar equation
\begin{multline}\label{mKdVUni}
\be^T(\D_\bk{\bf A}+\D_\bw{\bf B})\be V_T+\frac{1}{2}\be^T\bigg(\D_\bk^3{\bf B}(\be,\be,\be)-3\D_\bk^2{\bf B}(\be,{\bm \delta})\bigg)V^2V_X+\be^T{\bf K}V_{XXX} = 0\,.
\end{multline}
This completes the derivation of the mKdV equation.

\section{Application 1: Stratified Shallow Water}\label{sec:SWH}
The first example that the theory of this paper is applied to is a stratified shallow water system. Such a system forms a natural candidate for the theory, as it will be apparent that it possesses a doubly affine symmetry. This example will demonstrate how the solution associated with this symmetry, which turns out to be the uniform flow solution in each layer, can generate the required criticality for the mKdV to emerge. The mKdV has been derived in such settings in many works \cite{dr78,g99,gpt97,kb81} for the zero velocity background flow state, and so the theory of this paper allows one to take this one stage further to the case of finite background velocity states.

\begin{figure}
\includegraphics[width = 0.4\textwidth]{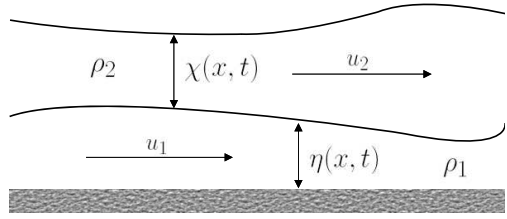}
\caption{A sketch of the system governed by the equations (\ref{SWSys}).}
\label{Fig:SWH}
\end{figure}

The shallow water model investigated here is rooted in the model proposed by \textsc{Baines} \cite{b88}, but is augmented with third order dispersive terms using the work of \textsc{Donaldson} \cite{d06}. This leads to the set of equations
\begin{equation}\label{SWSys}
\begin{split}
(\rho_1 \eta)_t+(\rho_1\eta u_1)_x&= 0\,,\\
(\rho _2 \chi)_t+(\rho _2\chi u_2)_x&=0\,,\\
(\rho _1 u_1)_t+\bigg(\frac{\rho_1}{2}u_1^2+g\rho_1\eta+g\rho_2\chi\bigg)_x &= a_{11}\eta_{xxx}+a_{12}\chi_{xxx}\,,\\
(\rho _2 u_2)_t+\bigg(\frac{\rho_2}{2}u_2^2 +g\rho_2\eta+g\rho_2\chi\bigg)_x  &= a_{21}\eta_{xxx}+a_{22}\chi_{xxx}\,.
\end{split}
\end{equation}
In the above $\rho_i,\,u_i$ denotes the density and velocity of the fluid in layer $i$, $g$ is acceleration due to gravity, $\eta$ denotes the thickness of layer 1, which is taken to be the lower layer and $\chi$ represents the thickness of the upper fluid in layer 2. For stable stratification, we impose that $\rho_2<\rho_1$, so that $r \equiv \frac{\rho_2}{\rho_1}<1$. The dispersive constants $a_{ij}$ are given by
\begin{equation}\notag
\begin{split}
a_{11} =& \sigma_1+\sigma_2-\frac{1}{3}\rho _1g\eta_0^2-\rho_2g\eta_0\chi_0-\frac{1}{2}g \chi _0^2,\\
a_{12} =& a_{21} = \sigma_2-\frac{1}{6}\rho_2 g \eta_0^2-\frac{1}{4}\rho_2g \eta_0 \chi _0-\frac{\rho_2^2}{2\rho_1}g \eta_0 \chi _0-\frac{5}{12}\rho_2 g \chi _0^2,\\
a_{22} =&\sigma_2 -\frac{\rho_2^2}{2 \rho_1}g \eta_0 \chi_0-\frac{1}{3}\rho_2 g \chi_0^2\,,
\end{split}
\end{equation}
where $\sigma_i$ denotes the surface tension constant for each fluid and the zero subscript denotes the quiescent thickness of the flow. This setup is pictured in figure \ref{Fig:SWH}

Under the assumption that the flow is irrotational, one may introduce the velocity potentials $\phi_i$ with the property that
\[
(\phi_i)_x = u_i\,,
\]
which allows one to then write (\ref{SWSys}) in potential form by integration of the last two equations with respect to $x$:
\begin{subequations}\label{SWSphi}
\begin{align}
(\rho_1 \eta)_t+(\rho_1\eta u_1)_x&= 0\,,\label{consSW1}\\
(\rho _2 \chi)_t+(\rho _2\chi u_2)_x&=0\,,\label{consSW2}\\
(\rho _1 \phi_1)_t+\frac{\rho_1}{2}\phi_1^2+g\rho_1\eta+g\rho_2\chi &=R_1+ a_{11}\eta_{xx}+a_{12}\chi_{xx}\,,\\
(\rho _2 \phi_2)_t+\frac{\rho_2}{2}\phi_2^2 +g\rho_2\eta+g\rho_2\chi  &= R_2+a_{21}\eta_{xx}+a_{22}\chi_{xx}\,.
\end{align}
\end{subequations}
This system is generated from the Lagrangian
\[
\begin{split}
\mathscr{L}  = & \iint \mathcal{L}(\phi_1,\phi_2,\eta,\chi,(\phi_1)_x,(\phi_2)_x,\eta_x,\chi_x) \, dx dt\\
 =& \iint \bigg[\rho_1\bigg(\eta(\phi_1)_t+\frac{1}{2}\eta(\phi_1)_x^2+\frac{g}{2}\eta^2\bigg)\\
 &+\frac{a_{11}}{2}\eta_x^2+a_{12}\eta_x\chi_x+\frac{a_{22}}{2}\chi_x^2-R_1\eta-R_2\chi\\
&\qquad+\rho_2\bigg(\chi(\phi_2)_t+\frac{1}{2}\chi(\phi_2)_x^2+g \eta \chi+\frac{g}{2}\chi^2\bigg) \bigg]\, dx\,dt\,.
\end{split}
\]

The system (\ref{SWSphi}) possesses two symmetries, one associated with the constant shift in each velocity potential. The solution associated with this symmetry is the constant flow solution in each layer, given by
\[
\phi_i = \theta_i\,.
\]
Substitution of this into (\ref{SWSys}) gives that the thicknesses for the uniform flow are given by
\begin{equation}\notag
\begin{split}
\eta_0 =& \frac{1}{g(\rho_1-\rho_2)}\bigg(\fr(\rho_2k_2^2-\rho_1k_1^2)+R_1-R_2-\rho_1\omega_1+\rho_2\omega_2\bigg),\\[2mm]
\chi_0 =& \frac{\rho_1}{g(\rho_1-\rho_2)}\bigg(R_2-R_1-\omega_2+\omega_1+\fr(k_1^2-k_2^2)\bigg)\,,
\end{split}
\end{equation}
where the $R_i$ result as constant of integration and can be thought of as Bernoulli constants for each layer. The above features, namely the generation of the problem from a Lagrangian density and the presence of a two parameter symmetry group, make the shallow water system (\ref{SWSphi}) a natural candidate to apply the theory of this paper.

\subsection{Conservation laws, criticality and the emergence of the mKdV}
\begin{figure}
\centering
\begin{subfigure}{0.45\textwidth}
\includegraphics[width=\textwidth,trim={0 0 1.5cm 0},clip]{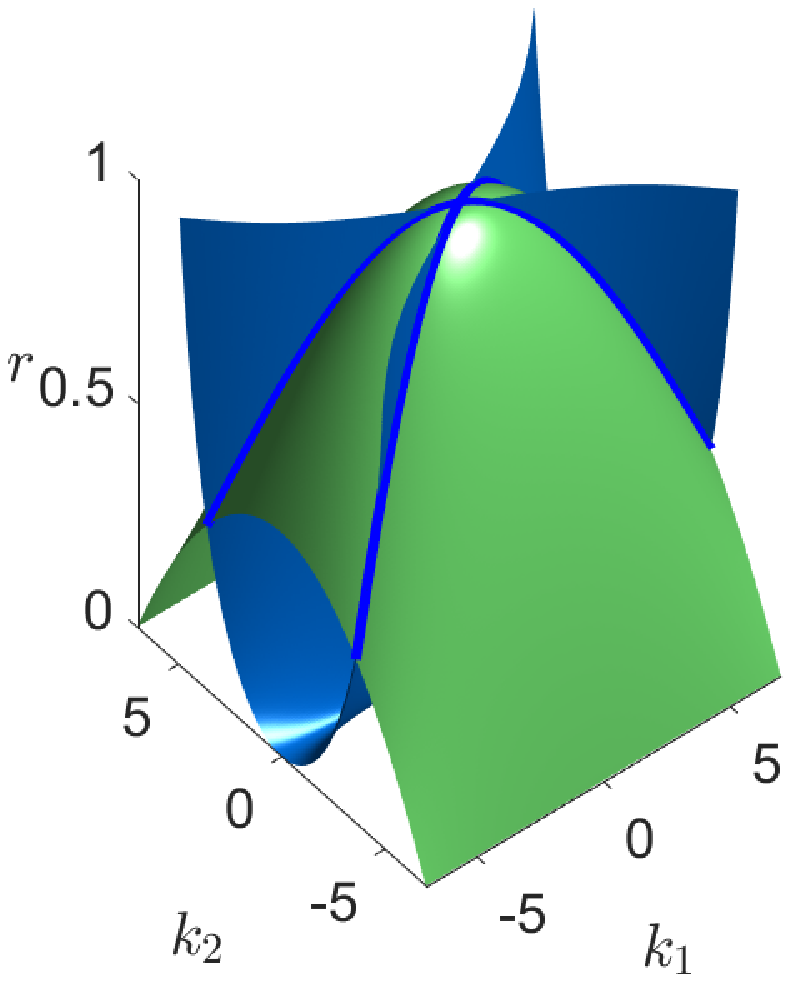}
\end{subfigure} \
\begin{subfigure}{0.45\textwidth}
\includegraphics[width=\textwidth,trim={0 0 0 0},clip]{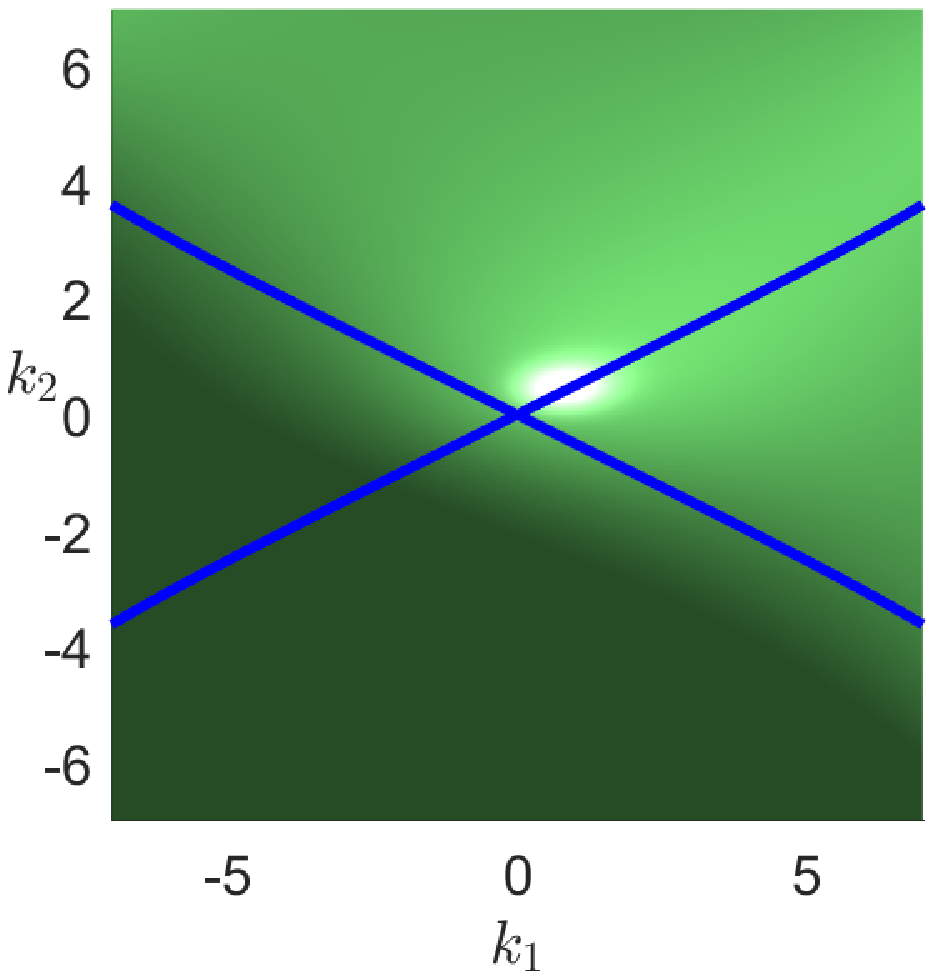}
\end{subfigure}
\caption{An illustration of how the criticality leading to the modified KdV may be met for $\eta_0 = 10,\,\chi_0 = 5$. The green surface indicates the surface where (\ref{SWHCrit1}) holds and the blue one represents (\ref{SWHCrit2}). Their intersection is highlighted with a blue line, with the modified KdV being the emergent modulation equation along it.}
\label{Fig:mKdVSWHCrit}
\end{figure}

The conservation laws for this system are given by (\ref{consSW1}) and (\ref{consSW2}), and so evaluated along the basic state the conservation law vectors are given by
\[
{\bf A} = 
\begin{pmatrix}
\rho_1\eta_0\\
\rho_2\chi_0
\end{pmatrix}\,, \quad 
{\bf B} = 
\begin{pmatrix}
\rho_1k_1\eta_0\\
\rho_2k_2 \chi_0
\end{pmatrix}\,.
\]
The first step in obtaining the mKdV for this system is to assess whether the relevant criticality conditions can be met. The first of these is met when
\begin{equation}\notag
\begin{split}
{\rm det}\big[ {\rm D}_{\rm k}{\bf B}\big] &= \ {\rm det}
\begin{pmatrix}
\rho_1\eta_0-\frac{\rho_1k_1^2}{g(1-r)} & \frac{\rho _2 k_1k_2}{g(1-r)}\\
\frac{\rho _2 k_1k_2}{g(1-r)} & \rho_2\chi_0-\frac{\rho_2k_2^2}{g(1-r)}
\end{pmatrix}
= 0\,,
\end{split}
\end{equation}
which can be reduced to
\begin{equation}\label{SWHCrit1}
(1-F_1^2)(1-F_2^2) = r\,.
\end{equation}
This expression arises from the literature of shallow water stratification as a stability boundary, but also corresponds to one of the system's characteristic speed vanishing \cite{b54,l90}. Providing this condition holds, it allows one to define the eigenvector of the zero eigenvalue of $\D_\bk{\bf B}$ as
\[
\be = \begin{pmatrix}
-\rho_2k_1k_2\\
g\rho_1\eta_0(1-r-F_1^2)
\end{pmatrix}\,.
\]
The second criticality arises from the expression
\[
\begin{split}
\be^T\D_\bk^2{\bf B}(\be,\be) = 3 g^2\rho_1^3\rho_2k_2\eta_0^2(1-r-F_1^2)\big[\chi_0r(1-F_2^2)F_1^2-\eta_0(1-F_1^2)^2F_2^2\big]\,.
\end{split}
\]
This only vanishes for physically relevant scenarios when the term in the square brackets is zero, meaning that
\begin{equation}\label{SWHCrit2}
\chi_0r(1-F_2^2)F_1^2 = \eta_0(1-F_1^2)^2F_2^2\,.
\end{equation}
This can be combined with the condition (\ref{SWHCrit1}) to give the $r$ independent condition
\begin{equation}\label{SWHTotalCond}
\chi_0(1-F_2^2)^2F_1^2 = \eta_0(1-F_1^2)F_2^2\,.
\end{equation}
The question remains as to whether both (\ref{SWHCrit1}) and (\ref{SWHCrit2}) can be met simultaneously, and to demonstrate that these can both be satisfied we appeal to geometric arguments. For each fixed $\eta_0,\,\chi_0$ the conditions (\ref{SWHCrit1}) and (\ref{SWHCrit2}) can be visualised in $(r,k_1,k_2)$-space, and an example of this is pictured in figure \ref{Fig:mKdVSWHCrit}. It would appear that the both conditions are satisfied for a continuum of values, and so the mKdV may be obtained. In such cases, the vector $\bdel$ exists and can be computed as
\[
\begin{split}
{\bm \delta} = \ \frac{\rho_2k_1}{g\eta_0(1-r-F_1^2)}\big(-3\rho_2k_1^2k_2^2-2g\rho_2k_2^2\eta_0(1-r-F_1^2)+g^2\rho_1\eta_0^2(1-r-F_1^2)^2\big)
\begin{pmatrix}
1\\
0
\end{pmatrix}\,.
\end{split}
\]
All that remains is to compute its coefficients.

For the coefficient of the time derivative term, one has that
\[
\begin{split}
\be^T
({\rm D}_{\bf k}{\bf A}+{\rm D}_\bw{\bf B})\be =&-2g^2\rho_1^2\rho_2\chi_0\eta_0^2(1-r-F_1^2)\bigg[ \frac{k_1}{g \eta_0}(1-F_2^2)+\frac{k_2}{g \chi_0}(1-F_1^2)\bigg]\,,
\end{split}
\]
which was also obtained elsewhere \cite{rb16a}. The coefficient of the dispersive term involves a Jordan chain argument, however the details of this appear elsewhere \cite{r17c} and simply state the result that
\[
\be^T {\bf K} = g\rho_1^2\eta_0^2\chi_0(1-r-F_1^2)\big(a_{11}r(1-F_2^2)-2ra_{12}+(1-F_1^2)a_{22}\big)\,.
\] 
The final component to compute is the coefficient of the cubic nonlinearity. This is done in two parts, with the first giving the result
\[
\begin{split}
\be^T\D_\bk^3{\bf B}(\be,\be,\be) =& \frac{3g^3\rho_1^5\eta_0^4(1-r-F_1^2)^4}{(1-r)} \big((1-F_1^2)(2r-1+F_1^2)-r\big)\,.
\end{split}
\]
The other term appearing in the cubic coefficient is given by
\[
\begin{split}
\be^T\D_\bk&{\bf B}(\bdel,\be)\\
 =& \ \frac{g^3\rho_1^3\rho_2^2\chi_0\eta_0^3(1-r-F_1^2)^2}{(1-r)}\big(2(1-r)+F_1^2F_2^2\big)^2\,.
\end{split}
\]
Combining these gives
\begin{equation}\notag
\begin{split}
\be^T\D_\bk^3{\bf B}(\be,\be,\be) -3\be^T\D_\bk^2{\bf B}(\bdel,\be)=& \ -3g^3\rho_2^2\rho_1^3\chi_0\eta_0^3(1-r-F_1^2)^2\big(F_1^2F_2^2+4(F_1^2+F_2^2)\big)\,.
\end{split}
\end{equation}
Therefore, by using the coefficient computed above one is able to construct the relevant mKdV as
\[
a_0 V_T+a_1V^2V_X+a_3V_{XXX} = 0\,,
\]
with
\[
\begin{array}{rcl}
a_0&=&\rho_2\bigg( \frac{k_1}{g \eta_0}(1-F_2^2)+\frac{k_2}{g \chi_0}(1-F_1^2)\bigg)\,,\\[2mm]
a_1&=&-\frac{3}{4}g\rho_1\rho_2 \eta_0F_2^2(1-F_1^2)\big(F_1^2F_2^2+4(F_1^2+F_2^2)\big)\,,\\[2mm]
a_2&=&-\frac{1}{2g}\big(a_{11}r(1-F_2^2)-2ra_{12}+(1-F_1^2)a_{22})\,.
\end{array}
\]
Noting that (\ref{SWHTotalCond}) implies that $1-F_1^2>0$, the sign of the nonlinear term appears to be in agreement with the zero velocity results \cite{dr78,kb81,g99}.

\section{Application 2: coupled nonlinear Schr\"odinger equations}\label{sec:NLS}
The second application presented, which presents a new emergence of the mKdV equation, is a set of coupled Nonlinear Schr\"odinger (NLS)
equations. Systems like this appear across a variety of contexts, such as when studying ocean waves \cite{ah15,lp10,o06,r76}, Bose-Einstein condensates \cite{ac13,sb09,yc12} and electromagnetic waves \cite{rp92}. Deriving nonlinear reductions like the mKdV in contexts such as the coupled NLS allows one to generate an analytic picture of the bifurcation of periodic travelling waves to various pairings of dark and bright solitary waves \cite{k90,r17b}, and so driving the mKdV in this context is of some interest.

The coupled NLS equations considered in this paper are given by
\begin{equation}\label{CNLS}
\begin{split}
i(\Psi_1)_t + \alpha _1 (\Psi_1)_{xx} + (\beta _{11}|\Psi_1|^2+\beta _{12}|\Psi_2|^2)\Psi_1 = 0\,,\\[3mm]
i(\Psi_2)_t + \alpha _2 (\Psi_2)_{xx} + (\beta _{21}|\Psi_1|^2+\beta _{22}|\Psi_2|^2)\Psi_2 = 0\,,
\end{split}
\end{equation}
for complex valued unknowns $\Psi_i(x,t)$ and $\alpha_i,\beta_{ij} \in \mathbb{R}$ constants. In order for this system to possess a generating Lagrangian density, we require $\beta_{12} = \beta_{21}$ and so in subsequent working we replace the latter with the former. In such a case, the Lagrangian which generates the set of equations (\ref{CNLS}) is given by
\[
\begin{split}
\mathscr{L} = \iint \frac{i}{2}&\big(\Psi_1^*(\Psi_1)_t-\Psi_1(\Psi_1)_t^*\big)+
 \frac{i}{2}\big(\Psi_2^*(\Psi_2)_t-\Psi_2(\Psi_2)_t^*\big)\\[2mm]
& -\alpha_1|(\Psi_1)_x|^2-\alpha_2|(\Psi_2)_x|^2
+\frac{1}{2}\beta_{11}|\Psi_1|^4 \\
&+ \beta_{12}|\Psi_1|^2|\Psi_2|^2 +
\frac{1}{2} \beta_{22}|\Psi_2|^4\ dx\,dt\,.
\end{split}
\]
The relative equilibrium solution is associated with the $SO(2)$ symmetries in each of the $\Psi_i$, which are independent. Associated with these are the plane wave solutions
\[
\Psi_i = \Psi_i^{(0)}e^{i \theta_i}\,,
\] 
and upon substitution into (\ref{CNLS}), one obtains that the amplitudes $\Psi_i^{(0)}$ satisfy
\[
\begin{split}
|\Psi_1^{(0)}|^2 &= \frac{1}{\beta}\big(\beta_{22}(\alpha_1k_1^2+\omega_1)-\beta_{12}(\alpha_2k_2^2+\omega_2)\big), \\ |\Psi_2^{(0)}|^2 &= \frac{1}{\beta}\big(\beta_{11}(\alpha_2k_2^2+\omega_2)-\beta_{12}(\alpha_1k_1^2+\omega_1)\big),
\end{split}
\]
where $\beta = \beta_{11}\beta_{22} - \beta_{12}^2$.

\subsection{Conservation laws, criticality and the emergent mKdV equation}
\begin{figure}
\centering
\begin{subfigure}{0.45\textwidth}
\includegraphics[width=\textwidth,trim={0 0 1.5cm 0},clip]{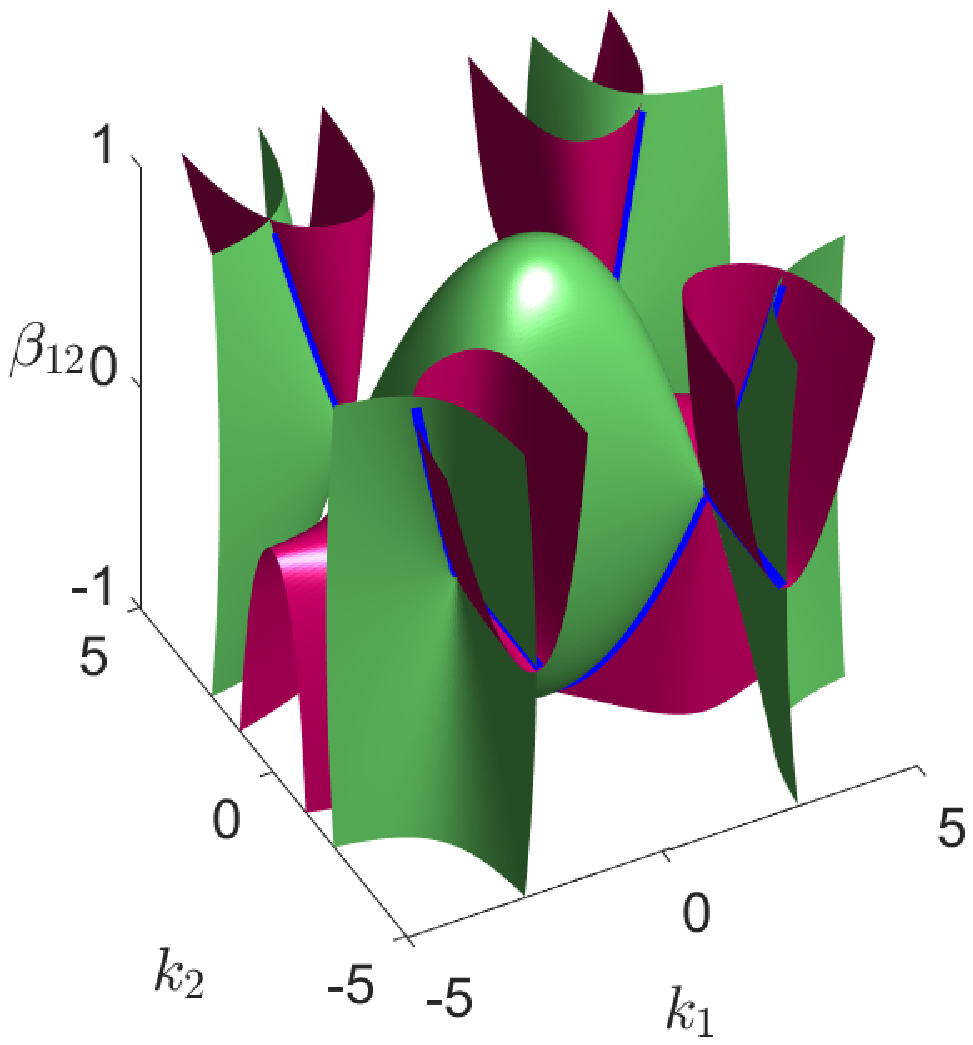}
\end{subfigure} \quad
\begin{subfigure}{0.45\textwidth}
\includegraphics[width=\textwidth,trim={0 0 0 0},clip]{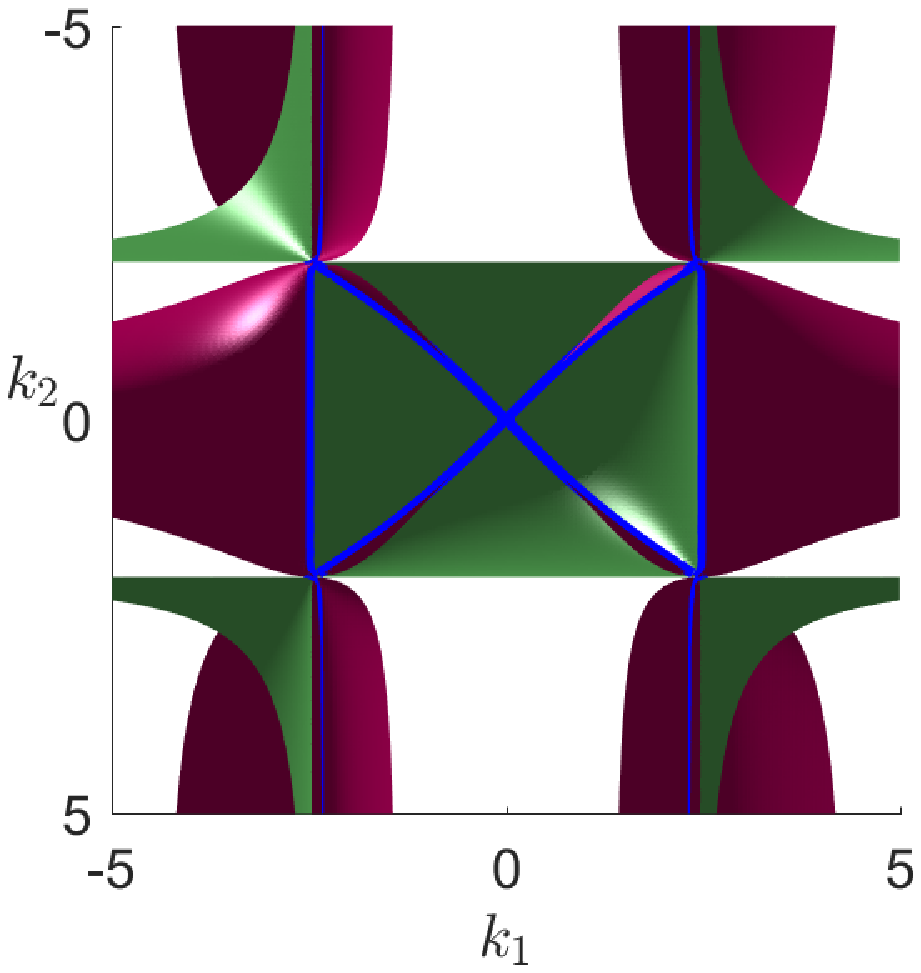}
\end{subfigure}
\caption{An illustration of how the criticality leading to the modified KdV may be met for $|\Psi_1^{(0)}| = 6,\,|\Psi_2^{(0)}| = 4,\,\beta_{11} = \beta_{22} = -1,\,\alpha_1 = \alpha_2 = \frac{1}{2}$. The green surface indicates the surface where (\ref{NLSCrit1}) holds and the blue one represents (\ref{NLSCrit2}). Their intersection is highlighted with a blue line, with the modified KdV being the emergent modulation equation along it.}
\label{Fig:mKdVNLSCrit}
\end{figure}

The conservation law components for the system (\ref{CNLS}) can be found as
\[
A = \frac{1}{2}\begin{pmatrix}
|\Psi_1|^2\\
|\Psi_2|^2
\end{pmatrix}\,, \quad B = \Im\begin{pmatrix}
(\Psi_1)_x\Psi_1^*\\
(\Psi_2)_x\Psi_2^*\\
\end{pmatrix}\,,
\]
where $*$ denotes the complex conjugate of the expression and $\Im$ denotes that the imaginary part of the expression is taken. We can evaluate these on the relative equilibrium solution to obtain the tensors required for the theory:
\begin{equation}
{\bf A} = \frac{1}{2}
\begin{pmatrix}
|\Psi_1^{(0)}|^2\\
|\Psi_2^{(0)}|^2
\end{pmatrix}\,, \quad {\bf B} = 
\begin{pmatrix}
k_1|\Psi_1^{(0)}|^2\\
k_2|\Psi_2^{(0)}|^2
\end{pmatrix}\,.
\end{equation}
These may be used to determine the relevant criticality required for the paper. The first occurs when the determinant of $\D_\bk{\bf B}$ vanishes, which explicitly means
\[
\D_\bk {\bf B} = \frac{1}{\beta}
\begin{pmatrix}
\alpha_1|\Psi_1^{(0)}|^2(1+\beta_{22}E_1^2)&-\frac{2\alpha_1\alpha_2k_1k_2\beta_{12}}{\beta}\\
-\frac{2\alpha_1\alpha_2k_1k_2\beta_{12}}{\beta}&\alpha_2|\Psi_2^{(0)}|^2(1+\beta_{22}E_2^2)
\end{pmatrix}\,,
\]
where to lighten the expressions we have introduced the dimensionless quantities
\[
E_1^2 = \frac{2\alpha_1k_1^2}{\beta|\Psi_1^{(0)}|^2}, \quad E_2^2 = \frac{2\alpha_2k_2^2}{\beta|\Psi_2^{(0)}|^2}\,.
\]
Simplification of this determinant leads to the expression
\begin{equation}\label{NLSCrit1}
(\beta_{11}+\beta E_1^2)(\beta_{22}+\beta E_2^2) = \beta_{12}^2\,.
\end{equation}
This forms the primary criticality condition, and has been shown to correspond to a stability boundary for the plane waves \cite{blp01}. The second criticality that must be met for the mKdV equation to emerge is
\begin{multline}\notag
\be^T{\bf B}(\be,\be) = \frac{6\alpha_1^3\alpha_2^2k_2|\Psi_1^{(0)}|^4(1+\beta_{22}E_1^2)}{\beta}\bigg(|\Psi_1^{(0)}|^2(1+\beta_{22}E_1^2)(\beta_{11}+\beta E_1^2)\\-\beta_{12}|\Psi_2^{(0)}|^2(1+\beta_{11}E_2^2)\bigg) = 0\,.
\end{multline}
This occurs when the term within the largest bracket vanishes. Therefore this condition requires that
\begin{equation}\label{NLSCrit2}
|\Psi_1^{(0)}|^2(1+\beta_{22}E_1^2)(\beta_{11}+\beta E_1^2)-\beta_{12}|\Psi_2^{(0)}|^2(1+\beta_{11}E_2^2) = 0,.
\end{equation}
A visualisation of when these coefficients are met simultaneously for fixed amplitudes is given in figure \ref{Fig:mKdVNLSCrit}.
When these are satisfied, the vector ${\bm \delta}$ can then be found to be
\begin{multline}\notag
{\bm \delta} = \frac{2\alpha_1^2\alpha_2 k_1}{\beta \beta_{12}}\bigg(2\beta_{12}|B_0|^2(\beta_{22}+\beta E_2^2)+\beta|A_0|^2(1+\beta_{22}E_1^2)^2\bigg)
\begin{pmatrix}
1\\
0
\end{pmatrix}\,.
\end{multline}

All that remains is to compute the relevant coefficients for the emerging mKdV equation. Starting with the coefficient of the time derivative, one has
\begin{multline}\notag
\be^T\big(\D_\bk {\bf A}+\D_\bw{\bf B})\be = \frac{2\alpha_1^2\alpha_2|\Psi_1^{(0)}|^2(1+\beta_{22}E_1^2)}{\beta}\bigg(|\Psi_2^{(0)}|^2(\beta_{22}+\beta E_2^2)k_1+|\Psi_1^{(0)}|^2(\beta_{11}+\beta E_1^2)k_2\bigg)\,.
\end{multline}
The next coefficient considered is that of the dispersive term. The full details of the Jordan chain analysis appear elsewhere \cite{r17b}, and lead to the result that
\begin{multline}\notag
\be^T{\bf K} = \frac{\alpha_1^2\alpha_2|\Psi_1^{(0)}|^2(1+\beta_{11}E_1^2)}{2\beta}\bigg(\alpha_2|\Psi_1^{(0)}|^2(\beta_{11}+\beta E_1^2)+\alpha_1|\Psi_2^{(0)}|^2(\beta_{22}+\beta E_2^2)\bigg)\,.
\end{multline}
Only the coefficient of the cubic nonlinearity remains to be computed. The first term considered is
\[
\be^T\D_\bk^3{\bf B}(\be,\be,\be) =  \frac{6\alpha_2^2\zeta_2^4}{\beta_{12}^2}\bigg(\beta E_1^2(1+\beta_{22}E_1^2)+(\beta_{11}+\beta E_1^2)\bigg)\,.
\]
The other component required for this coefficient is given by
\[
\begin{split}
\be^T\D_\bk^2{\bf B}(\bdel,\be)= \frac{2\alpha_1^4\alpha_2^2E_1^2|\Psi_1^{(0)}|^4(1+\beta_{22}E_1^2)}{\beta_{12}^2}\bigg(2\beta_{12}|\Psi_2^{(0)}|^2(\beta_{22}+\beta E_2^2)+|\Psi_1^{(0)}|^2\beta(1+\beta_{22}E_1^2)^2\bigg)^2\,.
\end{split}
\]
Combining these gives
\[
\begin{split}
\be^T\D_\bk^3{\bf B}(\be,\be,\be) -3\be^T\D_\bk^2{\bf B}(\bdel,\be) =& \ \frac{6\alpha_2^2\alpha_1^4\beta_{12}|\Psi_1^{(0)}|^6|\Psi_2^{(0)}|^2(1+\beta_{22}E_1^2)^2}{\beta}\big(3(\beta_{22}E_1^2+\beta_{11}E_2^2)-1\big)\,.\\
\end{split}
\]
Therefore, by using the coefficient computed above the modified KdV is given by
\[
a_0 V_T+a_1V^2V_X+a_3V_{XXX} = 0\,,
\]
with
\[
\begin{array}{rcl}
a_0 &=& |\Psi_2^{(0)}|^2(\beta_{22}+\beta E_2^2)k_1+|\Psi_1^{(0)}|^2(\beta_{11}+\beta E_1^2)k_2\,,\\[2mm]
a_1&=&\frac{3}{2}\alpha_2\alpha_1^2\beta_{12}|\Psi_1^{(0)}|^2|\Psi_2^{(0)}|^2(1+\beta_{22}E_1^2)\big(3(\beta_{22}E_1^2+\beta_{11}E_2^2)-1\big)\,,\\[2mm]
a_2 &=& \frac{1}{4}\bigg(\alpha_2|\Psi_1^{(0)}|^2(\beta_{11}+\beta E_1^2)+\alpha_1|\Psi_2^{(0)}|^2(\beta_{22}+\beta E_2^2)\bigg)\,.
\end{array}
\]

\section{Concluding remarks}
This paper has demonstrated that, if given a Lagrangian density whose Euler-Lagrange equations possess a two phase relative equilibria, the mKdV equation may be obtained providing suitable conditions are met. Moreover, an additional method to obtain the coefficient of the resulting nonlinearity was demonstrated and is in agreement with the one obtained from the reduction.

The multiphase analogy of the method of Kuramoto is a valuable step forward for multiphase modulation, since it allows one to deduce what the coefficients of the nonlinear terms are \emph{a priori}. This is beneficial in cases where the computation of the nonlinear coefficients becomes involved within the modulation. One expects this method to be invaluable in future analyses. For example, the method predicts that when the coefficient of the time derivative term in (\ref{mKdVUni}) vanishes, the relevant modulation equation should be
\begin{widetext}
\begin{multline}
\begin{split}
\be^T\big(\D_\bw {\bf A}\be-(\D_\bk{\bf A}+\D_\bw {\bf B})\bgam\big)V_{TT}+\bigg(\frac{1}{6}\be^T\big(\D_\bk^3{\bf B}(\be,\be,\be)-6\D_\bk^2{\bf B}(\be,{\bm \delta})\big)V^3+\be^T{\bf K}V_{XX}\bigg)_{XX}\\
+\be^T\bigg(\D_\bk^2{\bf A}(\be,\be)+\D_\bk\D_\bw{\bf B}(\be,\be)-\D_\bk^2{\bf B}(\be,\bgam)-\big(\D_\bk{\bf A}+\D_\bw{\bf B}){\bm \delta}\bigg)(VV_T)_X\\
+\be^T\bigg(\D_\bk\D_\bw{\bf B}(\be,\be)-\D_\bk{\bf B}(\be,\bgam)\bigg)(V_X\partial_X^{-1}V_T)_X = 0\,,
\end{split}
\end{multline}
\end{widetext}
where
\[
\D_\bk{\bf B}\bgam = (\D_\bk{\bf A}+\D_\bw{\bf B})\be\,,
\]
and $\partial_X^{-1}$ denotes the antiderivative. The above is a modified two-way Boussinesq, and the derivation of this via modulation will appear in another work.

The paper has only discussed the case of two symmetries, but the formulation of the problem allows this to be extended to arbitrarily many so long as the zero eigenvalue of $\D_\bk{\bf B}$ is simple. The case where the zero eigenvalue is nonsimple and the kernel of $\D_\bk{\bf B}$ has more than one element has the potential to lead to coupled nonlinear equations. This is because the projection of the final vector system can be done using each of these kernel elements. It remains to be answered whether one can generate a system of coupled mKdV equations, as well as the form in which these will emerge, and so further study is needed in this direction in order to answer this.

\section*{Acknowledgements}
The author would like to thank Prof. Tom Bridges for his discussions and insight during the formulation of this work. The author was in receipt of a fully funded Ph.D studentship under the EPSRC grant EP/L505092/1 during the formulation of this work.

\appendix

\section{Coefficient of the quadratic nonlinearity}\label{app:QuadCoeff}
Here we provide the details of how the coefficient of the quadratic nonlinearity is computed. This gives the result
\[
\begin{split}
\lth \Zh_{\theta_i},&\D^3S(\Zh)(\Zh_{k_j},\xi_5) - {\bf J}(\xi_5)_{\theta_j}
-\zeta_1{\bf J}\Zh_{k_1k_j} -\zeta_2{\bf J}\Zh_{k_jk_2}\rth\\[2mm]
=&\,\lth \D^3S(\Zh)(\Zh_{k_j},\Zh_{\theta_i})-{\bf J}\Zh_{\theta_i\theta_j},\xi_5\rth
\\&-\lth\Zh_{\theta_i},\zeta_1{\bf J}\Zh_{k_1k_j}+\zeta_2{\bf J}\Zh_{k_jk_2}\rth\,,\\[2mm]
=& -\lth\Zh_{\theta_ik_j},{\bf L}\xi_5\rth
-\lth\Zh_{\theta_i},\zeta_1{\bf J}\Zh_{k_1k_j}+\zeta_2{\bf J}\Zh_{k_jk_2}\rth\,,\\[2mm]
=& -\lth\Zh_{\theta_ik_j},\zeta_1{\bf J}\Zh_{k_1}+\zeta_2{\bf J}\Zh_{k_2}\rth
\\&-\lth\Zh_{\theta_i},\zeta_1{\bf J}\Zh_{k_1k_j}+\zeta_2{\bf J}\Zh_{k_jk_2}\rth\,,\\[2mm]
=&\,\zeta_1\partial_{k_1k_j}\mathscr{B}_i+\zeta_2\partial_{k_jk_2}\mathscr{B}_i.
\end{split}
\]
where we have used that
\[
{\bf L}\Zh_{\theta_i k_j} ={\bf J}\Zh_{\theta_i\theta_j} -\D^3S(\Zh)(\Zh_{k_j},\Zh_{\theta_i}) 
\,,
\] 
seen by differentiating (\ref{ThetaDeriv}) with respect to $k_j$.
Overall, this gives that the tensor acting on the nonlinearity takes the form 
\begin{multline}\notag
\begin{pmatrix}
\sum_{i=1}^2\zeta_i(\zeta_1\partial_{k_1k_i}\mathscr{B}_1 + \zeta_2\partial_{k_ik_2}\mathscr{B}_1) \\[2mm]
\sum_{i=1}^2\zeta_i(\zeta_1\partial_{k_1k_i}\mathscr{B}_2 + \zeta_2\partial_{k_ik_2}\mathscr{B}_2)
\end{pmatrix}\equiv\D^2_\bk{\bf B}(\be,\be)\,.
\end{multline}

\section{Details of the vanishing quadratic terms}\label{app:Vanishquad}
Here we provide the details leading to the zero coefficients of the quadratic terms in \S\ref{sec:ForuthOrd}.
Starting with the $U_XU_{XXX}$ term:
\begin{equation}\label{VanishQuad1}
\begin{split}
&\lth \Zh_{\theta_p},{\bf J}\kappa+\sum_{i=1}^2{\bf J}(\xi_6)_{\theta_i}-\D^3S(\Zh)(\xi_6,\Zh_{\theta_i})\rth\\
=&\sum_{i=1}^2\bigg[\zeta_i\bigg(-\lth\Zh_{k_p}, {\bf J}(\xi_5)_{\theta_i}-\D^3S(\Zh)(\xi_5,\Zh_{k_i})\\
&+\sum_{j=1}^2\zeta_j{\bf J}\Zh_{k_ik_j}\rth+\lth {\bf J}\xi_5,\Zh_{\theta_p k_i}\rth\bigg)+\delta_i\lth \Zh_{k_p},{\bf J}\Zh_{k_i}\rth\bigg]\,,\\
=& \sum_{i=1}^2\bigg[\zeta_i\big(-\lth\xi_5, {\bf L}\Zh_{k_p k_i}\rth-\sum_{j=1}^2\zeta_j\lth \Zh_{k_p},{\bf J}\Zh_{k_ik_j}\rth\big)\\
&+\delta_i\lth \Zh_{k_p},{\bf J}\Zh_{k_i}\rth\bigg]\,,\\
=& \sum_{i,j=1}^2\bigg[-\zeta_i\zeta_j\big(\lth{\bf J}\Zh_{k_j},\Zh_{k_p k_i}\rth+\lth \Zh_{k_p},{\bf J}\Zh_{k_ik_j}\rth\big)\\
&+\delta_i\lth \Zh_{k_p},{\bf J}\Zh_{k_i}\rth\bigg]\,,\\
=& \sum_{i,j=1}^2\bigg[-\zeta_i\zeta_j\partial_{k_i}\lth{\bf J}\Zh_{k_j},\Zh_{k_p}\rth+\delta_i\lth \Zh_{k_p},{\bf J}\Zh_{k_i}\rth\bigg] = 0\,.
\end{split}
\end{equation}
where we have used the result (\ref{k-k-zero}). Namely, this result highlights that
\[
\lth \Zh_{\theta_p},{\bf J}\kappa\rth = \sum_{i=1}^2\zeta_i\lth \xi_5,{\bf J}\Zh_{\theta_p k_i}\rth\,.
\]
This will be used in the computation of the coefficient of the $U_{XX}^2$ term:
\begin{equation}\label{VanishQuad2}
\begin{split}
&\lth \Zh_{\theta_p},{\bf J}\kappa-\frac{1}{2}\D^3S(\Zh)(\xi_5,\xi_5)\rth\\
  =& \frac{1}{2}\lth \xi_5,2\sum_{i=1}^2\zeta_i{\bf J}\Zh_{\theta_p k_i}-\D^3S(\Zh)(\xi_5,\Zh_{\theta_p})\rth\,,\\
=&\frac{1}{2}\lth \xi_5,\sum_{i=1}^2\zeta_i{\bf J}\Zh_{\theta_p k_i}\rth+\frac{1}{2}\lth \xi_5,{\bf L}(\xi_5)_{\theta_p}\rth\,,\\
=&\frac{1}{2}\sum_{i=1}^2\lth \xi_5,\zeta_i{\bf J}\Zh_{\theta_p k_i}\rth+\frac{1}{2}\lth \zeta_i{\bf J}\Zh_{k_i},(\xi_5)_{\theta_p}\rth\,,\\
 =& \frac{1}{2}(1-1)\sum_{i=1}^2\lth \xi_5,\zeta_i{\bf J}\Zh_{\theta_p k_i}\rth = 0\,.
\end{split}
\end{equation}
Thus, both terms which would be considered dissipative do not appear in the final PDE.

\section{Coefficient of the cubic nonlinearity}\label{app:CubicCoeff}
We provide the details of the calculation of the cubic coefficient of (\ref{mKdVUni}).  This will be done in stages, by first considering the terms in (\ref{fourthordmkdv}) containing $\kappa$:
\begin{widetext}
\[
\begin{split}
&\sum_{i=1}^2\zeta_i\lth \Zh_{\theta_p},{\bf J}(\kappa)_{\theta_i}-\D^3S(\Zh)(\kappa,\Zh_{k_i})\rth = \sum_{i=1}^2\zeta_i\lth\kappa,{\bf J}\Zh_{\theta_i\theta_p}-\D^3S(\Zh)(\Zh_{\theta_p},\Zh_{k_i})\rth = \sum_{i=1}^2\zeta_i\lth \kappa,{\bf L}\Zh_{\theta_p k_i}\rth\,,\\
=&\sum_{i=1}^2\zeta_i\lth \Zh_{\theta_p k_i},\sum_{j=1}^2\bigg[\zeta_j\big({\bf J}(\xi_5)_{\theta_j}-\D^3S(\Zh)(\Zh_{k_j},\xi_5)+\sum_{m=1}^2\zeta_m{\bf J}\Zh_{k_jk_m}\big)-\delta_j{\bf J}\Zh_{k_j}\bigg]\rth\,,\\
=&\sum_{i,j,m=1}^2\zeta_i\zeta_j\zeta_m\lth \Zh_{\theta_p k_i},{\bf J}\Zh_{k_jk_m}\rth-\sum_{i,j=1}^2\zeta_i\delta_j\lth \Zh_{\theta_p,k_i},{\bf J}\Zh_{k_j}\rth+\sum_{i,j=1}^2\zeta_i\zeta_j\lth \xi_5,{\bf J}\Zh_{\theta_p \theta_j k_i}-\D^3S(\Zh)(\Zh_{\theta_p k_i},\Zh_{k_j})\rth\,.
\end{split}
\]
Combine these with the terms involving $\xi_5$:
\[
\begin{split}
&\sum_{i=1}^2\bigg[\zeta_i\lth \Zh_{\theta_p},{\bf J}(\kappa)_{\theta_i}-\D^3S(\Zh)(\kappa,\Zh_{k_i})\rth-\sum_{i,j=1}^2\zeta_i\lth \Zh_{\theta_p},\frac{3}{2}\delta_j{\bf J}\Zh_{k_ik_j}\rth\\
&+\frac{1}{2}\zeta_j\D^3S(\Zh)(\xi_5,\Zh_{k_ik_j}) -\frac{1}{2}\D^4S(\Zh)(\xi_5,\Zh_{k_i},\Zh_{k_j})\rth-\frac{1}{2}\delta_i\lth \Zh_{\theta_p},{\bf J}(\xi_5)_{\theta_i}-\D^3S(\Zh)(\xi_5,\Zh_{k_j})\rth\bigg]\,,\\
=&\sum_{i,j,m=1}^2\zeta_i\zeta_j\zeta_m\lth \Zh_{\theta_p k_i},{\bf J}\Zh_{k_jk_m}\rth-\sum_{i,j=1}^2\zeta_i\delta_j\bigg(\lth \Zh_{\theta_p,k_i},{\bf J}\Zh_{k_j}\rth+\frac{3}{2}\lth \Zh_{\theta_p},\Zh_{k_ik_j}\rth\bigg)\\
&-\frac{1}{2}\sum_{i=1}^2\delta_i\lth \xi_5,{\bf L}\Zh_{\theta_p k_i}\rth+\sum_{i,j=1}^2\zeta_i\zeta_j\lth \xi_5,{\bf J}\Zh_{\theta_p \theta_j k_i}-\D^3S(\Zh)(\Zh_{\theta_p k_i},\Zh_{k_j})\\
&-\frac{1}{2}\D^3S(\Zh)(\Zh_{\theta_p},\Zh_{k_ik_j})-\frac{1}{2}\D^4S(\Zh)(\Zh_{\theta_p},\Zh_{k_i},\Zh_{k_j})\rth\,,\\
=&\sum_{i,j,m=1}^2\zeta_i\zeta_j\zeta_m\lth \Zh_{\theta_p k_i},{\bf J}\Zh_{k_jk_m}\rth-\frac{3}{2}\sum_{i,j=1}^2\zeta_i\delta_j(\lth \Zh_{\theta_p,k_i},{\bf J}\Zh_{k_j}\rth+\lth \Zh_{\theta_p},\Zh_{k_ik_j}\rth)+\frac{1}{2}\sum_{i,j=1}^2\zeta_i\zeta_j\lth {\bf L}\xi_5,\Zh_{\theta_p k_ik_j}\rth\,,\\
=&\frac{1}{2}\sum_{i,j,m=1}^2\zeta_i\zeta_j\zeta_m\big(2\lth \Zh_{\theta_p k_i},{\bf J}\Zh_{k_jk_m}\rth+\lth{\bf J}\Zh_{k_m},\Zh_{\theta_p k_ik_j}\rth\big)+\frac{3}{2}\sum_{i,j=1}^2\zeta_i\delta_j\partial_{k_i}\partial_{k_j}\mathscr{B}_p\,,\\
=&\frac{3}{2}\sum_{i,j,m=1}^2\zeta_i\zeta_j\zeta_m\lth \Zh_{\theta_p k_i},{\bf J}\Zh_{k_jk_m}\rth+\frac{3}{2}\sum_{i,j=1}^2\zeta_i\delta_j\partial_{k_i}\partial_{k_j}\mathscr{B}_p\,,\\
\end{split}
\]
where we have used the permutation of indices in the last step. Combination with the last term gives that
\[
\begin{split}
&\left\langle\!\left\langle \Zh_{\theta_p},\sum_{i=1}^2\Bigg[-\frac{1}{2}\delta_i\big(({\bf J}(\xi_5)_{\theta_i}-\D^3S(\Zh)(\xi_5,\Zh_{k_j})\big)+\zeta_i\bigg[{\bf J}(\kappa)_{\theta_i}-\D^3S(\Zh)(\kappa,\Zh_{k_i})\right.\right.\\
&-\sum_{j=1}^2\bigg(\delta_j{\bf J}\Zh_{k_ik_j}+\frac{1}{2}\zeta_j\D^3S(\Zh)(\xi_5,\Zh_{k_ik_j}) -\frac{1}{2}\D^4S(\Zh)(\xi_5,\Zh_{k_i},\Zh_{k_j})\left.\left.-\frac{1}{2}\sum_{m=1}^2\zeta_m\Zh_{k_i k_j k_m}\bigg)\bigg]\Bigg]\right\rangle\!\right\rangle\\
=&\frac{1}{2}\sum_{i,j,m=1}^2\zeta_i\zeta_j\zeta_m\big(3\lth \Zh_{\theta_p k_i},{\bf J}\Zh_{k_jk_m}\rth+\lth \Zh_{\theta_p},\Zh_{k_ik_jk_m}\rth\big)+\frac{3}{2}\sum_{i,j=1}^2\zeta_i\delta_j\partial_{k_i}\partial_{k_j}\mathscr{B}_p\,,\\
=&-\frac{1}{2}\sum_{i,j,m=1}^2\zeta_i\zeta_j\zeta_j\partial_{k_ik_jk_j}\mathscr{B}_p+\frac{3}{2}\sum_{i,j=1}^2\zeta_i\delta_j\partial_{k_ik_j}\mathscr{B}_p\,.
\end{split}
\]
\end{widetext}
Therefore,
\[
{\bf E} = \frac{1}{2}\big(3\D_\bk^2{\bf B}(\bdel,\be) -\D_\bk^3{\bf B}(\be,\be,\be)\big)\,,
\]
which matches the result obtained in (\ref{cubickura}).

\bibliographystyle{amsplain}

\end{document}